\documentclass[11pt]{article}
\usepackage{amssymb}
\usepackage[numbers]{natbib}
\usepackage{url}
\usepackage{graphicx}
\usepackage{xcolor}
\usepackage[figuresright]{rotating}
\usepackage{multirow}
\usepackage[margin=1in]{geometry}
\usepackage{setspace}
\usepackage[ruled,vlined,linesnumbered]{algorithm2e}

\SetKwComment{Comment}{$\triangleright$\ }{}
\SetSideCommentRight
\SetKw{Break}{break}

\usepackage{amsmath}
\allowdisplaybreaks[2]
\usepackage{amsthm}
\usepackage{amsfonts}
\usepackage{epic,eepic}
\usepackage{epsfig}
\usepackage{booktabs}
\usepackage[toc,page,title,titletoc,header]{appendix}
\usepackage{caption}
\usepackage{longtable}
\usepackage{pdflscape}
\usepackage{mathtools}
\mathtoolsset{showonlyrefs}
\usepackage{tabularx}
\usepackage{microtype}
\usepackage{subcaption}
\usepackage{placeins}
\usepackage{bm}
\usepackage{tikz}
\usepackage[hypertexnames=false]{hyperref}

\DeclareMathOperator*{\argmin}{arg\,min}
\DeclareMathOperator*{\argmax}{arg\,max}

\theoremstyle{plain}
\newtheorem{theorem}{Theorem}[section]
\newtheorem{proposition}[theorem]{Proposition}
\newtheorem{lemma}[theorem]{Lemma}

\newtheorem{definition}[theorem]{Definition}
\newtheorem{assumption}{Assumption}[section]

\newtheorem{condition}{Condition}[section]
\newtheorem{fact}[theorem]{Fact}
\theoremstyle{definition}
\newtheorem{setting}{Setting}[section]
\theoremstyle{plain}

\newcommand{\R}{\mathbb R}

\newcommand{\cO}{{\mathcal O}}
\newcommand{\tO}{\widetilde{\cO}}

\newcommand{\cL}{\mathcal L}
\newcommand{\APGz}{\hyperref[lem:conv-fistalinesearch]{\ensuremath{{\rm APG}_{0}}}}
\newcommand{\APGm}{\hyperref[lem:apgmu-complexity]{\ensuremath{{\rm APG}_{\mu}}}}
\newcommand{\IntV}{\hyperref[alg:inter-searchphi]{\ensuremath{{\rm IntV}}}}
\newcommand{\Bisec}{\hyperref[alg:bisec-multiphi]{\ensuremath{{\rm Bisec}}}}
\newcommand{\Dual}{\hyperref[alg:dual-algo]{\ensuremath{{\rm Dual}}}}
\newcommand{\VTest}{\hyperref[eq:dynamic-dual-call]{\ensuremath{{\rm VTest}}}}
\newcommand{\BiVFA}{\hyperref[alg:dynamic_bivfa]{\ensuremath{\rm BiVFA}}}
\newcommand{\SNtVFA}{\hyperref[alg:dynamic-newton]{\ensuremath{\rm SNtVFA}}}
\newcommand{\LBInit}{\hyperref[alg:lb-initialization]{\ensuremath{{\rm LB\_Init}}}}

\newcommand{\T}{\mathrm{T} }

\makeatletter

\newcommand{\Rmnum}[1]{\expandafter\@slowromancap\romannumeral #1@}
\makeatother

\title{Value-Function Root-Finding Algorithms for Composite Convex Simple Bilevel Optimization}
\author{Rujun Jiang\thanks{School of Data Science, Fudan University. Shanghai Key Laboratory for Contemporary Applied Mathematics, Fudan University. Emails: \texttt{rjjiang@fudan.edu.cn}} 
\and Xu Shi\thanks{School of Data Science, Fudan University. Email: \texttt{xshi22@m.fudan.edu.cn}} 
\and Weizheng Song\thanks{Department of Information Engineering, The Chinese University of Hong Kong. Email: \texttt{sw025@ie.cuhk.edu.hk}}
\and Jiulin Wang\thanks{School of Mathematical Sciences and LPMC, Nankai University. Email: \texttt{wangjiulin@nankai.edu.cn}}}
\date{}

\begin{document}

\onehalfspacing

\maketitle
\begin{abstract}
This paper considers composite convex simple bilevel optimization (SBO),
which minimizes a composite convex function over the solution set of another
composite convex minimization problem. We use a scalar value-function reformulation under which the bilevel optimal
value is characterized as the leftmost root of a scalar equation. To evaluate
this value function without level-set proximal access, we construct a
first-order oracle based on Lagrangian duality.
The
resulting feasibility and value-error bounds, together with multiplier
information, enable a bisection method and a safeguarded Newton-type method to
obtain an \((\epsilon,\epsilon)\)-solution. Under the stated assumptions and for
fixed problem data and initialization tolerance, both methods achieve an
operation complexity of \(\tO(\epsilon^{-1/2})\), matching, up to logarithmic
factors, the known near-optimal dependence for smooth convex SBO and the
accelerated first-order rate for unconstrained composite convex optimization.  
\end{abstract}

\section{Introduction}
\label{sec:introduction}

Simple bilevel optimization (SBO) concerns the problem of optimizing an
upper-level objective over the solution set of a lower-level problem
\citep{dempe2021simple}.
Such models arise in dictionary learning
\citep{cao2023projection,jiang2023conditional}, lifelong learning
\citep{malitsky2017chambolle,jiang2023conditional}, and safe or bi-objective
learning formulations \citep{kissel2020neural,gong2021bi}. Even in convex settings, the feasible set is given implicitly by lower-level
optimality, and projection onto it is generally unavailable.

In this paper, we consider the composite convex problem
\begin{equation}
\label{p:primal}\tag{SBO}
\min\limits_{\bm x\in \R^n} f(\bm x)\qquad{\rm s.t.}\quad \bm x\in \arg\min\limits_{\bm z\in\R^n} g(\bm z),
\end{equation}
where \(f=f_1+f_2\) and \(g=g_1+g_2\). The functions
\(f_1,g_1:\R^n\to\R\) are convex and continuously differentiable. 
Their gradients are Lipschitz continuous with constants \(L_{f_1}\) and
\(L_{g_1}\), respectively, assumed without loss of generality to satisfy
\(L_{f_1}>0,~ L_{g_1}\ge1\). The functions
\(f_2,g_2:\R^n\to\R\cup\{+\infty\}\) are proper, lower semicontinuous (l.s.c.),
and convex.
We assume that the lower-level and SBO solution sets are nonempty.

Let \(g^*:=\min_{\bm x\in{\rm dom}(g_2)}g(\bm x)\) and
\(f^*:=\min_{\bm x\in{\rm dom}(g_2)}f(\bm x)\), and let \(p^*\) denote the
optimal value of~\eqref{p:primal}. Given \(\epsilon>0\), we seek a one-sided
\((\epsilon,\epsilon)\)-solution, namely, a point
\(\widehat{\bm x}\in{\rm dom}(g_2)\) satisfying
\[
f(\widehat{\bm x})-p^*\le\epsilon
\qquad\text{and}\qquad
g(\widehat{\bm x})-g^*\le\epsilon.
\]

Following Wang et al.~\cite{wang2024near}, we consider the scalar value function
\begin{equation}
\label{eq:intro-scalar-value-function}
\bar g(c):=\min\{g(\bm x)\mid f(\bm x)\le c\}.
\end{equation}
The bilevel optimal value \(p^*\) is the leftmost root of
\(\bar g(c)=g^*\), reducing the computation of \(p^*\) to a one-dimensional scalar search.
Building on this characterization, we develop two value-function
search algorithms for the resulting scalar root-finding problem.
They access \(\bar g\) only through inexact oracle outputs
with checkable error bounds.
The bisection method \BiVFA{} maintains a search interval and tightens the
oracle tolerance whenever the test at the current scalar is inconclusive. The safeguarded Newton-type method \SNtVFA{} uses the returned multiplier to construct an affine minorant of
\(\bar g(c)-g^*\), whose zero defines a safeguarded scalar
step.

To implement both procedures, we construct a first-order value-function oracle
for the single-constraint composite subproblem associated with each trial
scalar \(c\). Each oracle call returns a primal point, checkable feasibility
and value-error bounds, and multiplier information that together support the
scalar decisions. A one-time initialization supplies a strict lower endpoint,
a reference point with a uniform strict-feasibility margin over all subsequent trial levels, and a uniform multiplier bound for all subsequent oracle
calls. These ingredients accommodate nonsmooth composite terms at both levels
without requiring upper-level strong convexity or level-set proximal access.
For fixed problem data and initialization tolerance, the initialization cost
is independent of the final accuracy, and each proposed algorithm computes an
\((\epsilon,\epsilon)\)-solution using \(\tO(\epsilon^{-1/2})\) basic
operations, as defined in Definition~\ref{def:basic-operation}.

\subsection{\texorpdfstring{Related Work}{Related Work}}

One line of work on convex SBO studies
regularization-based approximation schemes. Tikhonov-type, slowly vanishing, and
explicit hierarchical schemes replace SBO by
\(\min_{\bm x}\, g(\bm x)+\lambda_k f(\bm x)\), with \(\lambda_k>0\) vanishing
or chosen adaptively
\citep{Tikhonov1977,cabot2005proximal,solodov2007explicit}. Finite-time analyses obtain
\(\cO(\max\{\epsilon_f^{-1/(0.5-b)},\epsilon_g^{-1/b}\})\), \(b\in(0,0.5)\),
for related averaging iterative regularization schemes
\citep{kaushik2021method}, and
\(\cO(\max\{\epsilon_f^{-1/(1-a)},\epsilon_g^{-1/a}\})\), \(a\in(0.5,1)\),
for composite bi-subgradient schemes
\citep{merchav2023convex}. Additional structure sharpens
these rates: weak sharpness yields \(\tO(\epsilon^{-1/2})\)-type dependence
\citep{samadi2024achieving}, while H\"olderian
error bounds give explicit polynomial rates
\citep{chen2024penalty}.
For composite objectives at both levels, FBi-PG (dynamic FISTA) establishes
simultaneous lower- and upper-level objective rates, with an accelerated
lower-level rate and sharper H\"olderian-error-bound guarantees
\citep{merchav2026dynamic}, while iterative-regularization analyses provide
analogous simultaneous rates and extensions to accelerated schemes and surrogate
formulations \citep{shtern2026convergence}.
Further diagonal analyses cover proximal-gradient iterations for
composite SBO \citep{latafat2025convergence} and continuous-time-guided schemes
that yield last-iterate rates and weak convergence under Attouch--Czarnecki or
H\"olderian-growth conditions \citep{bot2026accelerating}.
These approximation schemes use assumptions, oracle models, or
asymptotic regimes different from the proximal-composite oracle setting with
checkable oracle bounds considered here.

A second line develops several other first-order SBO schemes under various
structural assumptions and oracle models. Minimal-norm and sequential averaging
methods handle smooth or composite variants under upper-level
smoothness and strong-convexity assumptions
\citep{beck2014first,sabach2017first}. Classical bundle
and \(\epsilon\)-subgradient methods address nonsmooth convex SBO through
cutting-plane or subgradient operations
\citep{solodov2007bundle,helou2017}. ITALEX-type methods handle
nonsmooth, non-strongly convex outer objectives by combining level-set expansion
with first-order inner solves \citep{doron2023methodology}, while related
first-order schemes based on online convex optimization \citep{shen2023online},
conditional-gradient or projection-free updates
\citep{jiang2023conditional,cao2023projection,giang2024projection}, and
accelerated gradients for smooth SBO \citep{cao2024accelerated} use
projection-type access, linear-optimization oracles, or cutting-plane
constructions. The present work complements these studies by handling nonsmooth terms through composite proximal operations and providing finite-time guarantees for both upper-level optimality and lower-level value gaps.

The most closely related methods implement bisection over scalar
optimal-value functions. For a trial upper-level value \(c\), \citep{wang2024near} considers
\(\bar g(c)\) in~\eqref{eq:intro-scalar-value-function} and uses that \(p^*\)
is the leftmost root of \(\bar g(c)=g^*\). This scalar value-function bisection approach
was first proposed in that work and attains near-optimal complexity by approximately solving
\(\min_{\bm x}\, g_1(\bm x)+g_2(\bm x)+{\rm I}_{{\rm Lev}_f(c)}(\bm x)\),
where \({\rm Lev}_f(c):=\{\bm x\mid f(\bm x)\le c\}\), and, for
any set \(\mathcal C\), \({\rm I}_{\mathcal C}\) denotes the indicator
function that equals \(0\) on \(\mathcal C\) and \(+\infty\) otherwise.
This formulation is natural when proximal access to
\(g_2+{\rm I}_{{\rm Lev}_f(c)}\) is available. In contrast, our
value-function oracle does not require this level-set proximal access.
It solves Lagrangian-dual subproblems using the combined proximal
operation \(\operatorname{prox}_{\tau(g_2+\theta f_2)}\), rather than the
proximal mapping of \(g_2+{\rm I}_{{\rm Lev}_f(c)}\).
Subsequent work on functionally constrained problems develops related scalar root-finding methods and lower bounds
in smooth and Lipschitz settings
\citep{zhang2024functionally}. Thus, the closest
\(\tO(\epsilon^{-1/2})\)-type comparisons differ in problem class and oracle
model: \citep{wang2024near} handles composite SBO through a
level-set proximal oracle, whereas \citep{zhang2024functionally} studies smooth
or Lipschitz functionally constrained formulations. 

A related but separate issue is the upper-level criterion.
Because near-feasible points may satisfy \((f(\bm x)<p^*,\, g(\bm x)>g^*)\), the one-sided condition
\(f(\bm x)-p^*\le\epsilon\) is weaker than the absolute, or two-sided,
condition \({ \lvert f(\bm x)-p^*\rvert}\le\epsilon\). This distinction is made explicit in
\citep{zhang2024functionally}, which proves general intractability of absolute  
solutions, while
\citep{giang2025convergence} characterizes asymptotic duality and
well-posedness conditions under one-sided guarantees without a finite-time
complexity bound.

\paragraph{Organization.}

Section~\ref{sec:preliminaries} gives the problem setting, the
value-function reformulation, and the black-box guarantees used by the
proposed algorithms. Sections~\ref{sec:sbo-target}
and~\ref{sec:newton-algorithm}
develop the bisection and safeguarded Newton-type
value-function search methods, respectively,
including their correctness and operation-complexity analyses. Section~\ref{sec:oracle-construction} develops a generic dual oracle and specializes it to construct the value-function oracle. Section~\ref{sec:auxiliary} constructs the reference certificate that
initializes both methods and analyzes its one-time cost.
Section~\ref{sec:experiment} reports numerical experiments,
and Section~\ref{sec:conclusion} concludes. The appendix contains
the detailed proof of the generic dual-oracle guarantee.

\paragraph{Notation.}
We use \({\rm I}_{\mathcal C}\), \({\rm dist}(\cdot,\mathcal C)\),
\({\rm dom}(h)\), and \(\partial h\) in their standard convex-analysis senses.
The norm \(\|\cdot\|\) is Euclidean for vectors and spectral for matrices. For
\(a\in\R\), let \([a]_+:=\max\{a,0\}\), and let \(\lceil a\rceil_+\) be the
smallest nonnegative integer no smaller than \(a\).

\section{\texorpdfstring{Preliminaries and Oracle Guarantees}{Preliminaries and Oracle Guarantees}}
\label{sec:preliminaries}

This section sets up the value-function framework used throughout the paper.
We state the standing assumption, derive the scalar root-finding reformulation
of SBO, and collect the subroutine guarantees on which the
proposed algorithms rely. These guarantees include the corresponding operation-complexity
bounds.

\subsection{\texorpdfstring{Standing Assumption and Value-Function Reformulation}{Standing Assumption and Value-Function Reformulation}}
\label{subsec:setting-value-reformulation}

For the value-function analysis, we impose the following
compactness, compatibility, and nondegeneracy conditions.

\begin{assumption}
\label{ass:basic}
\noindent\textup{(i)} The functions \(f_2\) and \(g_2\) are
proper, l.s.c., and convex, and the restriction of \(g_2\) to
\({\rm dom}(g_2)\) is continuous.
The set \({\rm dom}(g_2)\) is nonempty, compact,
and convex. Moreover, \({\rm dom}(g_2)\subseteq{\rm dom}(f_2)\), and \(f_2\)
is \(l_{f_2}\)-Lipschitz continuous on \({\rm dom}(g_2)\).

\noindent\textup{(ii)} The optimal value \(p^*\)
of Problem~\eqref{p:primal} satisfies \(p^*>f^*\).
\end{assumption}

Assumption~\ref{ass:basic} is natural for constrained SBO models. The compact
convex set \({\rm dom}(g_2)\) is the lower-level domain, and
\({\rm dom}(g_2)\subseteq{\rm dom}(f_2)\) ensures that the upper-level
nonsmooth term is finite on all relevant points, so \(f_2\) needs to be
Lipschitz continuous only there. The continuity condition on \(g_2\) is also
relative to \({\rm dom}(g_2)\), so it covers indicator terms and finite
continuous lower-level regularizers. If boundedness is imposed through \(f_2\)
instead, the lower-level domain may be replaced by
\({\rm dom}(g_2)\cap{\rm dom}(f_2)\), discarding only points with upper-level
value \(+\infty\). The condition \(p^*>f^*\) excludes the degenerate case in
which a minimizer of \(f\) over \({\rm dom}(g_2)\) is already lower-level
optimal. 

Under Assumption~\ref{ass:basic}, the lower-level and SBO problems both attain
their optimal values. We fix the following diameter and upper-variation bounds:
\begin{equation}
\label{eq:Dg-Omegag}
D_g\ge {\rm diam}({\rm dom}(g_2)),\qquad
\Omega_g:=\sup_{\bm x\in{\rm dom}(g_2)}\{g(\bm x)-g^*+1\}<\infty.
\end{equation}
The continuity of
\(g_1\) and of \(g_2\) on \({\rm dom}(g_2)\), together with the compactness of
\({\rm dom}(g_2)\), ensures that \(\Omega_g<\infty\). The quantity
\(\Omega_g\) is used only in complexity estimates.

To state operation-complexity bounds, we use the following notion of a basic operation.
\begin{definition}[Basic operation]
\label{def:basic-operation}
A \emph{basic operation} is one evaluation of \(f_1, f_2, g_1, g_2\),
\(\nabla f_1\), \(\nabla g_1\),
\(\operatorname{prox}_{\tau(g_2+\theta f_2)}\), or
\(\operatorname{prox}_{\tau(f_2+{\rm I}_{{\rm dom}(g_2)})}\), where
\(\tau>0\) and \(\theta\ge0\).
\end{definition}
In this definition, the parameter \(\theta\) is the multiplier parameter in the value-function subproblems,
and we assume access to \(\operatorname{prox}_{\tau(g_2+\theta f_2)}\), the combined nonsmooth
proximal call associated with weighted lower- and upper-level terms, as in
regularization and penalty methods for SBO
\citep{solodov2007explicit,merchav2023convex,shtern2026convergence,
merchav2026dynamic,samadi2024achieving,chen2024penalty,bot2026accelerating,latafat2025convergence}.
The latter proximal call is used for upper-level APG steps restricted to
\({\rm dom}(g_2)\).

We next recall the scalar value function \(\bar g\)
from~\eqref{eq:intro-scalar-value-function}, which is the basis of the
value-function search algorithms. The following facts, recalled from
\cite[Section 3.1]{wang2024near}, justify the scalar reformulation and the
bisection and Newton-type updates developed below.

\begin{fact}
\label{fact}
\noindent\textup{(i)} The function \(\bar g\) is convex and nonincreasing on
\((f^*,+\infty)\) \citep{Rockafellar1970convex}.

\noindent\textup{(ii)} If \(f^*<c<p^*\), then \(\bar g(c)>g^*\). If \(c\ge p^*\), then \(\bar g(c)=g^*\).

\noindent\textup{(iii)} The value \(p^*\) is the leftmost root of
\(\bar g(c)=g^*\).
\end{fact}

Fact~\ref{fact} identifies \(p^*\) as the leftmost zero of the nonincreasing
residual \(\Phi(c):=\bar g(c)-g^*\). This reduction suggests bisection and
Newton-type scalar searches, but the proposed methods obtain residual values and affine
information only through inexact value-function solves. At a trial level \(c\),
the relevant constrained composite subproblem is
\begin{equation}
\label{p:subp-rewrite}\tag{\({\rm VF}_c\)}
\min\limits_{\bm x\in \R^n} \ g(\bm x)=g_1(\bm x)+g_2(\bm x)
\qquad {\rm s.t.}\quad f_c(\bm x)\triangleq f(\bm x)-c\le0.
\end{equation}
An approximate primal solution of~\eqref{p:subp-rewrite} alone is not sufficient
for scalar search. The oracle must also provide feasibility and value
estimates for the residual comparisons in the scalar search, as well as
multiplier information for the Newton-type affine model.

\subsection{\texorpdfstring{Lower-Level APG Comparison Guarantee}{Lower-Level APG Comparison Guarantee}}
\label{subsec:oracle-interface}

Before the scalar search, the proposed algorithms obtain a lower-level comparison point with an
objective-gap guarantee using an accelerated proximal-gradient (APG) call
\citep{nesterov1983method,beck2009fast}. Subsequent
bisection and Newton-type routines use only this gap guarantee. The next lemma
records the corresponding guarantee, obtained by specializing
Lemma~\ref{lem:conv-fistalinesearch} to \(\varphi=g\).

\begin{lemma}[Lower-level APG guarantee]
\label{lem:front-apgz-certificate}

The routine \(\APGz\) is the non-strongly convex APG solver for the
lower-level composite problem
\(\min_{\bm x\in\R^n}\ g(\bm x)=g_1(\bm x)+g_2(\bm x)\).
Suppose that Assumption~\ref{ass:basic} holds, with \(D_g\)
fixed as in~\eqref{eq:Dg-Omegag}. For any
\(\delta_g>0\) and \(\bm x_{\rm init}^g\in{\rm dom}(g_2)\), the call
\(\bm x_g
=\APGz(g_1,g_2,\delta_g,\bm x_{\rm init}^g,D_g)\)
returns \(\bm x_g\in{\rm dom}(g_2)\) satisfying
\(0\le g(\bm x_g)-g^*\le \delta_g\),
and requires \(\cO(D_g\sqrt{L_{g_1}/\delta_g})\) basic operations.
\end{lemma}

\subsection{\texorpdfstring{Reference Certificate and Initialization Guarantee}{Reference Certificate and Initialization Guarantee}}
\label{subsec:front-lbinit}

Before any value-function test is made, the value-function search algorithms need a lower
endpoint for the scalar search, a reference point for certifying strict feasibility, and a strict-feasibility margin that is uniform
over subsequent trial levels. We collect these quantities in a reference certificate.
This subsection states the initialization guarantee used in
Sections~\ref{sec:sbo-target} and~\ref{sec:newton-algorithm}, while
Section~\ref{sec:auxiliary} gives the construction.

\begin{condition}[Reference certificate]
\label{cond:reference}
 With \(D_g\) and \(\Omega_g\) fixed as in~\eqref{eq:Dg-Omegag}, we call
the tuple \((l_0,\Delta,\bm x_f,M_\Delta)\) a reference certificate if
\[
\bm x_f\in{\rm dom}(g_2),\qquad f^*<l_0<p^*,\qquad
l_0=f(\bm x_f)+\Delta,\qquad \Delta>0,
\]
and there exists \(\bm x_g^{\rm ref}\in{\rm dom}(g_2)\) such that
\(g(\bm x_g^{\rm ref})-g^*\le1\) and
\begin{equation}
\label{eq:Mz-def}
\begin{aligned}
M_\Delta
\triangleq
\max\left\{
\frac{g(\bm x_f)-g(\bm x_g^{\rm ref})+1}{\Delta},\,1
\right\} 
\le
\max\left\{\frac{\Omega_g}{\Delta},\,1\right\},
\end{aligned}
\end{equation}
where the inequality follows from \(g(\bm x_g^{\rm ref})\ge g^*\) and the definition of $\Omega_g$.
\end{condition}

 For the initialization guarantee and value-function oracle, define
\begin{equation}
\label{eq:CD-def}
C_D:=1+D_g(1+D_g).
\end{equation}
Let $\delta_{\rm init}>0$
denote the prescribed tolerance for the initialization procedure. Using the midpoint between \(f^*\) and \(p^*\), define
\begin{equation}
\label{eq:m0-def}
m_0:=\bar g\bigl((p^*+f^*)/2\bigr)-g^*,
\qquad
Q_{\rm init}:=
\max\left\{\frac1{\delta_{\rm init}},\frac8{p^*-f^*},\frac8{m_0}\right\}.
\end{equation}
By Fact~\ref{fact} and Assumption~\ref{ass:basic}, the strict inequality \(p^*>f^*\) implies \(m_0>0\).
Together with \(\delta_{\rm init}>0\), this shows that \(Q_{\rm init}\) is well defined and positive.
We remark that the constant \(Q_{\rm init}\), which is used only in the complexity bounds, is not an input to the initialization procedure.
The rates below describe the dependence on \(\epsilon\) for fixed
\(p^*-f^*\), \(m_0\), and \(\delta_{\rm init}\).

The proposed algorithms also require a Lipschitz constant
\(B_f \ge1\) for \(f\) on
the lower-level domain. Under Assumption~\ref{ass:basic}, such a constant exists
by the compactness of \({\rm dom}(g_2)\), the continuity of
\(\nabla f_1\), and the Lipschitz continuity of \(f_2\) on
\({\rm dom}(g_2)\). The following lemma records this fact.
\begin{lemma}[Lipschitz continuity of \(f\)]
\label{lem:lip-bound-f}
Suppose that Assumption~\ref{ass:basic} holds. Then there exists a constant
\(B_f{\ge1}\) such that \(f\) is \(B_f\)-Lipschitz continuous on
\({\rm dom}(g_2)\).
\end{lemma}
Henceforth, \(B_f\) denotes any such known Lipschitz constant satisfying
\(B_f\ge1\). The next lemma records the guarantee
and cost for constructing a reference certificate.

\begin{lemma}[Initialization guarantee]
\label{lem:front-lbinit}

There exists an initialization procedure \(\LBInit\) such that, whenever
Assumption~\ref{ass:basic} holds, \(B_f \ge1\), \(f\) is
\(B_f\)-Lipschitz continuous on
\({\rm dom}(g_2)\), \(\bm x_{\rm init}^f,\bm x_{\rm init}^g\in{\rm dom}(g_2)\),
and \(0<\delta_{\rm init}\le{1/B_f}\), the call
\(\LBInit(\delta_{\rm init},D_g,B_f,\bm x_{\rm init}^f,\bm x_{\rm init}^g)\)
returns \((l_0,\Delta,\bm x_f,M_\Delta)\), a reference certificate satisfying
Condition~\ref{cond:reference}, with the estimates

\begin{equation}
\label{eq:MDelta-init-bound}
M_\Delta\le \Omega_g Q_{\rm init},
\qquad
\Delta\ge Q_{\rm init}^{-1}.
\end{equation}
 Its one-time initialization cost is
\[
T_{\rm pre}
=\cO\!\left(
\left[\sqrt{C_D L_{f_1}\Omega_g}\,Q_{\rm init}
+\sqrt{C_D L_{g_1}}\sqrt{Q_{\rm init}}\right]
{\left(1+\log Q_{\rm init}\right)^2}
\right).
\]
\end{lemma}

Section~\ref{sec:auxiliary} constructs \(\LBInit\) and proves Theorem~\ref{thm:lb-total}, which establishes the guarantee in Lemma~\ref{lem:front-lbinit}. The subsequent analysis uses
only Condition~\ref{cond:reference} and the bounds in
Lemma~\ref{lem:front-lbinit}.

\subsection{\texorpdfstring{Value-Function Oracle Guarantee}{Value-Function Oracle Guarantee}}
\label{subsec:value-oracle-interface}

We now state the value-function oracle guarantee used by the proposed algorithms.
The construction is deferred to Section~\ref{sec:oracle-construction}. The
bisection and Newton-type analyses use only the lemma below.

Fix a reference certificate \((l_0,\Delta,\bm x_f,M_\Delta)\) satisfying
Condition~\ref{cond:reference}. For \(f_c(\bm x):=f(\bm x)-c\) and every
\(c\ge l_0\), Condition~\ref{cond:reference} gives
\(f_c(\bm x_f)\le-\Delta\). Thus \(\bm x_f\) is strictly feasible for
\(f_c(\bm x)\le0\), and this margin underlies the multiplier bound \(M_\Delta\).
The oracle uses \(B_f\) to restrict the requested tolerance and \(C_D\)
from~\eqref{eq:CD-def} to set the perturbation scale. Its dual-oracle scaling
constant is
\begin{equation}
\label{eq:CDDelta}
C_\Delta
\triangleq
4C_D M_\Delta.
\end{equation}

For \(0<\eta\le{1/B_f}\), define the strongly convex perturbation

\begin{equation}
\label{eq:G-eta-def}
G_{\eta}(\bm x):=g(\bm x)+\frac{\eta}{2}\|\bm x-\bm x_f\|^2.
\end{equation}
The associated perturbed value-function problem is
\[
\min\limits_{\bm x\in\R^n}\ G_\eta(\bm x)
\qquad {\rm s.t.}\quad f_c(\bm x)\le0.
\]

The value-function test provides computable estimates of the scalar residual
\(\Phi(c)=\bar g(c)-g^*\). A direct evaluation of \(\bar g(c)\) would require
solving the constrained value problem~\eqref{p:subp-rewrite}. We instead use
the strict-feasibility margin from the reference certificate and the strongly
convex perturbation \(G_\eta\) in~\eqref{eq:G-eta-def} to invoke the generic
dual oracle in the basic-operation model of
Definition~\ref{def:basic-operation}. The resulting primal point yields approximate feasibility and two-sided control of \(\bar g(c)\). The associated
multiplier yields the dual lower bound used to form affine lower
estimates. For \(c\ge l_0\), \(0<\delta\le{1/B_f}\), and
\(\bm x_{\rm in}\in{\rm dom}(g_2)\), set \(\eta:=\delta/C_D\).
The value-function test depends on the fixed reference certificate through \(\bm x_f\),
\(M_\Delta\), and \(C_\Delta\). We suppress this dependence in the notation
and define the value-function test by the specialized dual-oracle call
\begin{equation}
\label{eq:dynamic-dual-call}
(\bm x_c,z_c)=\VTest(c,\delta,\bm x_{\rm in})
:=\Dual\left(G_\eta,f_c,1,M_\Delta,\bm x_{\rm in};\delta/C_\Delta\right).
\end{equation}

The next lemma records the resulting value estimate, dual lower bound,
and operation cost. Its proof is given in
Section~\ref{sec:pre-dual}.

\begin{lemma}[\(\VTest\) oracle guarantee]
\label{lem:front-dual-oracle}

Suppose that Assumption~\ref{ass:basic} holds. Let
\((l_0,\Delta,\bm x_f,M_\Delta)\) be a reference certificate satisfying
Condition~\ref{cond:reference}.
For every \(c\ge l_0\), \(0<\delta\le{1/B_f}\), and
\(\bm x_{\rm in}\in{\rm dom}(g_2)\), the call
\((\bm x_c,z_c)=\VTest(c,\delta,\bm x_{\rm in})\)
in~\eqref{eq:dynamic-dual-call} satisfies

\begin{equation}
\label{eq:front-dual-value-cert}
\bm x_c\in{\rm dom}(g_2),\quad
z_c\in[0,2M_\Delta),\quad
f_c(\bm x_c)\le\delta,\quad
{\lvert g(\bm x_c)-\bar g(c)\rvert}\le\delta.
\end{equation}
For the unperturbed dual value
\(d_c(z):=\inf_{\bm x\in{\rm dom}(g_2)}\{g(\bm x)+z f_c(\bm x)\}\),
the returned pair also satisfies
\begin{equation}
\label{eq:front-dual-lower}
\bar g(c)\ge d_c(z_c)>g(\bm x_c)-(D_g^2+2)\delta.
\end{equation}
The oracle call has operation complexity
\begin{equation}
\label{eq:front-dual-cost}
\cO\!\left(\sqrt{\frac{C_D(L_{g_1}+M_\Delta L_{f_1})}{\delta}}\,
{ \left(1+\log\frac1\delta\right)^2}\right).
\end{equation}
\end{lemma}

\section{Bisection Value-Function Algorithm for SBO}
\label{sec:sbo-target}

This section presents \BiVFA{}, the bisection value-function algorithm for
SBO. It first constructs a reference certificate through
\(\LBInit\), computes a lower-level comparison point by APG, and then performs a checked bisection search using
\(\VTest\).

\subsection{\texorpdfstring{Decision Rule and Bisection Algorithm}{Decision Rule and Bisection Algorithm}}
\label{subsec:bivfa-result-state}

The bisection phase uses \(\VTest\) through the following decision rule. A
comparison gap larger than \(\delta\) justifies a lower update, while a gap at
  most \(\epsilon/2\) yields an upper endpoint.

\begin{lemma}[Bisection decision guarantee]
\label{lem:front-dual-decision}

Suppose that Assumption~\ref{ass:basic} holds. Let
\((l_0,\Delta,\bm x_f,M_\Delta)\) be a reference certificate satisfying
Condition~\ref{cond:reference}. For any
\(c\ge l_0\), \(0<\delta\le{1/B_f}\), \(\epsilon>0\), and
\(\bm x_{\rm in}\in{\rm dom}(g_2)\), let
\((\bm x_c,z_c)=\VTest(c,\delta,\bm x_{\rm in})\) be defined
by~\eqref{eq:dynamic-dual-call}.
If \(\bm x_g\in{\rm dom}(g_2)\) satisfies
\(g(\bm x_g)-g^*\le\epsilon/2\), then the following decision implications hold:
\begin{align}
g(\bm x_c)-g(\bm x_g)>\delta
&\quad\Longrightarrow\quad c<p^*, \label{eq:front-dual-decision-lower}\\
g(\bm x_c)-g(\bm x_g)\le\epsilon/2
&\quad\Longrightarrow\quad
g(\bm x_c)\le g^*+\epsilon. \label{eq:front-dual-decision-upper}
\end{align}
\end{lemma}

\begin{proof}

Condition~\ref{cond:reference} gives \(l_0=f(\bm x_f)+\Delta\), so
\(c\ge l_0\) yields \(c\ge f(\bm x_f)+\Delta\). Hence
Lemma~\ref{lem:front-dual-oracle} applies to the call
in~\eqref{eq:dynamic-dual-call}, and the value estimate~\eqref{eq:front-dual-value-cert} holds.

To prove~\eqref{eq:front-dual-decision-lower}, suppose
\(g(\bm x_c)-g(\bm x_g)>\delta\). The value estimate \eqref{eq:front-dual-value-cert} gives
$
\bar g(c)\ge g(\bm x_c)-\delta>g(\bm x_g)\ge g^*,
$
and Fact~\ref{fact} gives \(c<p^*\). This proves
\eqref{eq:front-dual-decision-lower}. 
To prove
\eqref{eq:front-dual-decision-upper}, suppose
\(g(\bm x_c)-g(\bm x_g)\le\epsilon/2\). Combining this inequality
with \(g(\bm x_g)-g^*\le\epsilon/2\) yields
$
g(\bm x_c)
\le g(\bm x_g)+\epsilon/2
\le g^*+\epsilon$.
Thus~\eqref{eq:front-dual-decision-upper} follows.
\end{proof}

We now turn Lemma~\ref{lem:front-dual-decision} into a checked bisection
scheme for the residual equation \(\Phi(c)=\bar g(c)-g^*\). At iteration \(t\),
Algorithm~\ref{alg:dynamic_bivfa} (\BiVFA{}) maintains
\begin{equation}
\label{eq:bivfa-invariant}
l_0\le l_t<p^*,\qquad u_t=f(\bm x_u),\qquad
\bm x_u\in{\rm dom}(g_2),\qquad
g(\bm x_u)\le g^*+\epsilon,
\end{equation}
where \(l_t\) is a valid lower endpoint and \(\bm x_u\) gives the
current candidate output level \(u_t=f(\bm x_u)\).
The initialization computes a lower-level comparison point \(\bm x_g\) and sets
\(u_0=f(\bm x_g)\). If
\begin{equation}
\label{eq:bivfa-initial-certificate}
l_0<p^*,\qquad u_0=f(\bm x_g),\qquad u_0-l_0\le\epsilon,
\qquad \bm x_g\in{\rm dom}(g_2),\qquad
g(\bm x_g)\le g^*+\epsilon/2,
\end{equation}
then the displayed relations imply
\(f(\bm x_g)-p^*\le u_0-l_0\le\epsilon\) and
\(g(\bm x_g)-g^*\le\epsilon/2\). So the algorithm returns \(\bm x_g\) immediately as
an \((\epsilon,\epsilon)\)-optimal point. Otherwise, \(\VTest\) is called at each midpoint, and the decision rule either justifies a lower update, supplies a new
{upper endpoint}, or leaves an inconclusive comparison to be retested at a
smaller tolerance.

\begin{algorithm}[ht!]
\caption[Bisection value-function algorithm (\BiVFA{})]{Bisection value-function algorithm: \(\widehat{\bm x} = \BiVFA(\epsilon,D_g,B_f,
\delta_{\rm init},\bm x_{\rm init}^f,\bm x_{\rm init}^g)\)}
\label{alg:dynamic_bivfa}
\DontPrintSemicolon
\KwIn{\(D_g\) fixed as in~\eqref{eq:Dg-Omegag},
\(B_f \ge1\) a known Lipschitz constant for \(f\) on \({\rm dom}(g_2)\),
\(0<\epsilon\le\min\{1,4\delta_{\rm init}\}\),
\(0<\delta_{\rm init}\le{1/B_f}\),
\(\bm x_{\rm init}^f,\bm x_{\rm init}^g\in{\rm dom}(g_2)\).}
\KwOut{\(\widehat{\bm x}\) with
\(f(\widehat{\bm x})-p^*\le\epsilon\) and
\(g(\widehat{\bm x})-g^*\le\epsilon\).}
Invoke
\((l_0,\Delta,\bm x_f,M_\Delta)\gets
\LBInit(\delta_{\rm init},D_g,B_f,\bm x_{\rm init}^f,\bm x_{\rm init}^g)\).\label{alg:bivfa-lbinit}\;

Set \(C_D\) as in~\eqref{eq:CD-def} and \(C_\Delta\) as
in~\eqref{eq:CDDelta}.\;
Invoke \(\bm x_g{\gets}
\APGz
(g_1,g_2,\epsilon/2,\bm x_{\rm init}^g,D_g)\) and set
\(u_0\gets f(\bm x_g)\).\label{alg:dyn-main-uppercall}\label{alg:dyn-main-u0}\;
\If{\(u_0-l_0\le \epsilon\)}{
\Return \(\widehat{\bm x}=\bm x_g\).\label{alg:dyn-main-early}
\Comment*[r]{early objective-gap return}
}
Set \(t\gets0\), \(\bm x_u\gets\bm x_g\), and
\(\delta_0^{\rm in}\gets\delta_{\rm init}\).\label{alg:phase-init}\;
\While{\(u_t-l_t>\epsilon\)\label{alg:phase-outer}}{
Set \(w_t\gets u_t-l_t\), \(c_t\gets(l_t+u_t)/2\),
\(\delta_t^0\gets\min\{\delta_t^{\rm in},w_t/4\}\), and
\(j\gets0\).\label{alg:phase-midpoint}\label{alg:phase-epsinit}\;
Invoke
\((\bm x_{c_t},\cdot)\gets
\VTest(c_t,\delta_t^j,\bm x_u)\).\label{alg:phase-firstcall}\;
\While{\(\epsilon/2<g(\bm x_{c_t})-g(\bm x_g)\le\delta_t^j\)\label{alg:phase-inner}}{
\(j\gets j+1\) and \(\delta_t^j\gets\delta_t^{j-1}/2\).\label{alg:phase-refine}
\Comment*[r]{inconclusive comparison}
Invoke
\((\bm x_{c_t},\cdot)\gets
\VTest(c_t,\delta_t^j,\bm x_{c_t})\).\label{alg:phase-refinecall}\;
}
Set \(\delta_t\gets\delta_t^j\) and
\(\delta_{t+1}^{\rm in}\gets\delta_t\).\; \label{alg:phase-finaltol}
\uIf{\(g(\bm x_{c_t})-g(\bm x_g)>\delta_t\)}{
\((l_{t+1},u_{t+1})\gets(c_t,u_t)\).\label{alg:phase-lc}
\Comment*[r]{lower endpoint update}
}
\Else{
\(\bm x_u\gets\bm x_{c_t}\) and
\((l_{t+1},u_{t+1})\gets(l_t,f(\bm x_u))\).\label{alg:phase-uc}
\Comment*[r]{upper endpoint update}
}
\(t\gets t+1\).\; \label{alg:phase-next}
}
\Return \(\widehat{\bm x}=\bm x_u\).\label{alg:dyn-main-return}
\end{algorithm}

The next theorem gives the main result for \BiVFA{}. Its
proof is given in Subsection~\ref{subsec:dynamic-bisection}.
\begin{theorem}[\BiVFA{} guarantee]
\label{thm:correctness-bisection}
Suppose that Assumption~\ref{ass:basic} holds. Let \(D_g\) be fixed as
in~\eqref{eq:Dg-Omegag}, and let \(B_f \ge1\) be a Lipschitz constant for \(f\) on
\({\rm dom}(g_2)\). Then \BiVFA{} returns an
\((\epsilon,\epsilon)\)-optimal solution of Problem~\eqref{p:primal}, and, with
\(Q_{\rm init}\) and \(T_{\rm pre}\) defined in \eqref{eq:m0-def} and
Lemma~\ref{lem:front-lbinit}, respectively, its total operation complexity satisfies
\begin{equation}\label{eq:bisectionbd}
T_{\rm BiVFA}
=
T_{\rm pre}
+
\cO\!\left(
D_g\sqrt{\frac{L_{g_1}}{\epsilon}}
+
\sqrt{\frac{C_D}{\epsilon}
\left(L_{g_1}+L_{f_1}\Omega_g Q_{\rm init}\right)}
{ \left(1+\log\frac1\epsilon\right)^3}
\right).
\end{equation}

\end{theorem}

For fixed problem data and a fixed coarse preprocessing tolerance
\(\delta_{\rm init}\), \(T_{\rm pre}\), \(Q_{\rm init}\),
and all other problem-dependent factors hidden in the \(\cO\)-notation are independent
of the final accuracy \(\epsilon\). They can therefore be treated as constants
when interpreting the \(\epsilon\)-dependence of the bound.
Consequently,
Theorem~\ref{thm:correctness-bisection} gives
\(T_{\rm BiVFA}=\tO(\epsilon^{-1/2})\). Up to logarithmic factors,
this is the same accelerated first-order dependence as in unconstrained
composite convex optimization
\citep{nemirovskij1983problem,woodworth2016tight} and matches, up to logarithmic factors, the
\(\epsilon\)-dependence of the complexity lower bound for smooth convex SBO
established in~\citep{zhang2024functionally}. The bound also matches the rate of
recent level-set-proximal bisection methods for constrained value-function
formulations \citep{wang2024near,zhang2024functionally}. 
Our contribution is therefore on the oracle side: we retain the
\(\tO(\epsilon^{-1/2})\) dependence known for smooth SBO
\citep{zhang2024functionally} and level-set-proximal bisection
\citep{wang2024near}, but for nonsmooth composite \(f_2\) and \(g_2\), without
upper-level strong convexity \citep{beck2014first,sabach2017first} or proximal
access to \(g_2+{\rm I}_{{\rm Lev}_f(c)}\) \citep{wang2024near}.

\subsection{\texorpdfstring{Bisection-Phase Analysis}{Bisection-Phase Analysis}}
\label{subsec:dynamic-bisection}

Line~\ref{alg:bivfa-lbinit} calls \(\LBInit\) to construct the reference
certificate \((l_0,\Delta,\bm x_f,M_\Delta)\). We prove Theorem~\ref{thm:correctness-bisection} by conditioning
on this tuple. The bisection-phase analysis then reduces to two
estimates: contraction of the search interval and
the total number of tolerance refinements, which together control the number and
cost of \(\VTest\) calls. Condition~\ref{cond:reference} ensures that every such
call has \(c\ge l_0\) and
\(0<\delta\le\delta_{\rm init}\le{1/B_f}\), so the value estimates in~\eqref{eq:front-dual-value-cert} apply throughout. At each executed outer
iteration \(t\) (including \(t=0\) when the early return is not taken), set
\(w_t:=u_t-l_t\), \(c_t:=(l_t+u_t)/2\), and
\(\delta_t^0:=\min\{\delta_t^{\rm in},w_t/4\}\).

Refinement occurs only in the inconclusive region. The inner loop at
Line~\ref{alg:phase-inner} is entered precisely when
\(\epsilon/2<g(\bm x_{c_t})-g(\bm x_g)\le\delta_t^j\). In this case,
Line~\ref{alg:phase-refine} halves the tolerance and
Line~\ref{alg:phase-refinecall} repeats the midpoint test. Once the comparison
leaves this inconclusive region, Lemma~\ref{lem:front-dual-decision} gives a valid endpoint update:
\(g(\bm x_{c_t})-g(\bm x_g)>\delta_t\) justifies a lower update, while
\(g(\bm x_{c_t})-g(\bm x_g)\le\epsilon/2\) yields an {upper endpoint}.
The next lemma records the resulting interval contraction and the uniform lower
bound on all testing tolerances.

\begin{lemma}
\label{lem:gap-contract}

Suppose that \(0<\epsilon\le4\delta_{\rm init}\), and let \(N_{\rm out}\)
denote the number of outer iterations executed by
Algorithm~\ref{alg:dynamic_bivfa}. Then, at every outer iteration before
termination, all \(\VTest\) calls use testing tolerance at least
\(\epsilon/4\), and \(w_{t+1}\le 3w_t/4\). If \(u_0-l_0>\epsilon\), then
$
N_{\rm out}\le
\left\lceil\frac{\log_2((u_0-l_0)/\epsilon)}{\log_2(4/3)}\right\rceil_+.
$
\end{lemma}

\begin{proof}

We prove the tolerance claim by induction. Since
\(0<\epsilon\le4\delta_{\rm init}\), \(\delta_0^{\rm in}\ge\epsilon/4\).
At any executed outer iteration \(t\), the while condition in
Line~\ref{alg:phase-outer} gives \(w_t=u_t-l_t>\epsilon\). Hence, if
\(\delta_t^{\rm in}\ge\epsilon/4\), then
\(\delta_t^0=\min\{\delta_t^{\rm in},w_t/4\}\ge\epsilon/4\). A refinement can
occur only when
\(\delta_t^j\ge g(\bm x_{c_t})-g(\bm x_g)>\epsilon/2\), so the next tested tolerance satisfies
\(\delta_t^{j+1}=\delta_t^j/2>\epsilon/4\). Hence all calls in iteration \(t\), and the next input tolerance
\(\delta_{t+1}^{\rm in}=\delta_t\), are at least \(\epsilon/4\).

For contraction, \(\delta_t\le\delta_t^0\le w_t/4\). A lower update gives
\(w_{t+1}=w_t/2\). On an upper update, Line~\ref{alg:phase-uc} gives
\(u_{t+1}=f(\bm x_{c_t})\) and \(l_{t+1}=l_t\). Moreover,
applying \eqref{eq:front-dual-value-cert} to the call returning
\(\bm x_{c_t}\) gives
\[
f_{c_t}(\bm x_{c_t})\le\delta_t,\qquad
f(\bm x_{c_t})\le c_t+\delta_t\le l_t+\frac{3w_t}{4},
\]
where the last inequality uses \(c_t=l_t+w_t/2\) and
\(\delta_t\le w_t/4\). Thus
\(w_{t+1}\le3w_t/4\) in either branch. Iterating gives
\(w_t\le(3/4)^t(u_0-l_0)\), and the claimed \(N_{\rm out}\) bound follows.
\end{proof}

The uniform testing-tolerance bound in Lemma~\ref{lem:gap-contract} also
controls the number of tolerance refinements. The next lemma records this
refinement count.

\begin{lemma}
\label{lem:total-refinements}

Suppose that \(0<\epsilon\le4\delta_{\rm init}\) and \(u_0-l_0>\epsilon\).
Let \(J_{\rm ref}^{\rm Bi}\) be the total number of executions of
Line~\ref{alg:phase-refine}. Then
\(J_{\rm ref}^{\rm Bi}\le
\left\lceil\log_2(2\delta_{\rm init}/\epsilon)\right\rceil_+\).

\end{lemma}

\begin{proof}

Each execution of Line~\ref{alg:phase-refine} replaces
\(\delta_t^j\) by \(\delta_t^{j+1}=\delta_t^j/2\). Across outer iterations,
Lines~\ref{alg:phase-finaltol} and~\ref{alg:phase-epsinit} give
\(\delta_{t+1}^{\rm in}=\delta_t\) and
\(\delta_{t+1}^0=\min\{\delta_{t+1}^{\rm in},w_{t+1}/4\}\le\delta_t\).
Thus transitions between outer iterations never increase the testing tolerance, and the number of
refinements is maximized by a single halving sequence starting from
\(\delta_{\rm init}\). A refinement occurs only if
\(\epsilon/2<g(\bm x_{c_t})-g(\bm x_g)\le\delta_t^j\). So every counted
halving starts from a tolerance \(\delta_t^j>\epsilon/2\).
The number of such halvings is at most the claimed
\(J_{\rm ref}^{\rm Bi}\) bound.
\end{proof}

We complete the proof of
Theorem~\ref{thm:correctness-bisection} by verifying correctness of the
returned point, including early termination, and deriving the operation bound
from the value-function oracle cost.

\begin{proof}[Proof of Theorem~\ref{thm:correctness-bisection}]

\emph{Correctness.}
By Lemma~\ref{lem:front-lbinit}, Line~\ref{alg:bivfa-lbinit}
returns a tuple satisfying Condition~\ref{cond:reference}, with
preprocessing cost \(T_{\rm pre}\) and
\(M_\Delta\le\Omega_g Q_{\rm init}\). We condition on this tuple
throughout the bisection-phase analysis.

If Line~\ref{alg:dyn-main-early} is taken,
Condition~\ref{cond:reference} and Lemma~\ref{lem:front-apgz-certificate},
applied at Line~\ref{alg:dyn-main-uppercall}, give the relations in
\eqref{eq:bivfa-initial-certificate}, which is exactly the early-termination
condition above. Hence the returned point is \((\epsilon,\epsilon)\)-optimal.

Suppose now that the bisection loop is entered. We maintain
\eqref{eq:bivfa-invariant} throughout the loop. At \(t=0\), it follows from
Condition~\ref{cond:reference}, the assignments
\(u_0=f(\bm x_g)\) and \(\bm x_u=\bm x_g\), and
Lemma~\ref{lem:front-apgz-certificate} applied at
Line~\ref{alg:dyn-main-uppercall}. For the induction step, the tolerance
recursion keeps every tested tolerance in
\(0<\delta_t^j\le\delta_{\rm init}\le{1/B_f}\). The invariant gives
\(\bm x_u\in{\rm dom}(g_2)\) for the first call at each midpoint, and
Lemma~\ref{lem:front-dual-oracle} ensures that each refinement point remains in
\({\rm dom}(g_2)\). Together with
\(c_t>l_t\ge l_0>f^*\), these facts ensure that
Lemma~\ref{lem:front-dual-decision} applies to each midpoint call
with \(\bm x_c=\bm x_{c_t}\). In the
lower-update branch, the lemma gives
\[
g(\bm x_{c_t})-g(\bm x_g)>\delta_t=\delta_t^j
\quad\Longrightarrow\quad c_t<p^*,
\]
so \(l_{t+1}=c_t\) satisfies
\(l_0\le l_{t+1}<p^*\). On an upper update, the lower-update test fails and the
refinement loop exits through the upper-update branch. Since
\(\delta_t=\delta_t^j\), this gives
\[
g(\bm x_{c_t})-g(\bm x_g)\le\epsilon/2
\quad\Longrightarrow\quad g(\bm x_{c_t})\le g^*+\epsilon.
\]

Applying \eqref{eq:front-dual-value-cert} to the call returning
\(\bm x_{c_t}\) also gives
\(\bm x_{c_t}\in{\rm dom}(g_2)\).
Together with the assignment in Line~\ref{alg:phase-uc}, this verifies the
upper-endpoint part of
\eqref{eq:bivfa-invariant}.
Since the loop is entered only when \(u_0-l_0>\epsilon\), the
contraction estimate in Lemma~\ref{lem:gap-contract} implies finite
termination. 
{With \(N_{\rm out}\) defined in Lemma~\ref{lem:gap-contract}, the outer loop terminates at index \(t=N_{\rm out}\).}
At termination, \eqref{eq:bivfa-invariant} and the stopping
condition give 
\[f(\bm x_u)-p^*\le { u_{N_{\rm out}}-l_{N_{\rm out}}}\le\epsilon,\qquad g(\bm x_u)-g^*\le\epsilon.\]

\emph{Operation complexity.}
By Lemma~\ref{lem:front-apgz-certificate}, the APG call that
computes \(\bm x_g\) costs \(\cO(D_g\sqrt{L_{g_1}/\epsilon})\).
Let \(N_{\rm vf}^{\rm Bi}\) denote the total number of \(\VTest\) calls in the
bisection phase. Each outer iteration
performs one initial \(\VTest\) evaluation, and each refinement adds one call, so
\(N_{\rm vf}^{\rm Bi}\le N_{\rm out}+J_{\rm ref}^{\rm Bi}\).
Lemma~\ref{lem:total-refinements} gives
\({ J_{\rm ref}^{\rm Bi}=\cO(\log(1/\epsilon))}\) for fixed
\(\delta_{\rm init}\).

The initial width appearing in \(N_{\rm out}\) is controlled by the diameter
bound and the Lipschitz constant \(B_f\) for \(f\):
\[
u_0-l_0=f(\bm x_g)-f(\bm x_f)-\Delta
\le { \lvert f(\bm x_g)-f(\bm x_f)\rvert}\le B_fD_g.
\]
By Lemma~\ref{lem:gap-contract}, this gives
\({ N_{\rm out}=\cO(\log(1/\epsilon))}\) for fixed problem data.
Together with the call-count bound above and the bound on
\(J_{\rm ref}^{\rm Bi}\), this yields
\({ N_{\rm vf}^{\rm Bi}=\cO(\log(1/\epsilon))}\).
For every bisection-phase call
\(\VTest(c_t,\delta_t^j,\cdot)\), Lemma~\ref{lem:gap-contract} gives
\(\epsilon/4\le\delta_t^j\le\delta_{\rm init}\). Applying
\eqref{eq:front-dual-cost} with \(\delta=\delta_t^j\) therefore gives,
uniformly over all such calls,
\[
T_{\rm Dual}:=
\cO\!\left(
\sqrt{\frac{C_D(L_{g_1}+M_\Delta L_{f_1})}{\epsilon}}
{ \left(1+\log\frac1\epsilon\right)^2}
\right).
\]

Consequently, the bisection-phase cost is bounded by
\(\cO(D_g\sqrt{L_{g_1}/\epsilon})+N_{\rm vf}^{\rm Bi}T_{\rm Dual}\).
Substituting \(M_\Delta\le\Omega_g Q_{\rm init}\)
from~\eqref{eq:MDelta-init-bound} and adding the preprocessing cost
\(T_{\rm pre}\) gives \eqref{eq:bisectionbd}.
\end{proof}

\section{\texorpdfstring{Newton-Type Value-Function Algorithm for SBO}{Newton-Type Value-Function Algorithm for SBO}}
\label{sec:newton-algorithm}
This section develops \(\SNtVFA{}\), a safeguarded Newton-type value-function
algorithm for SBO. The method repeatedly calls \(\VTest\) to obtain certified primal--dual information, constructs affine minorants of the residual
\(\Phi(c)=\bar g(c)-g^*\), and uses these minorants to perform a safeguarded update of the trial level.

\subsection{\texorpdfstring{Safeguarded Newton-Type Search}{Safeguarded Newton-Type Search}}
\label{subsec:dynamic-newton}

We next introduce \(\SNtVFA{}\), summarized in
Algorithm~\ref{alg:dynamic-newton}, by deriving its scalar search step. By Fact~\ref{fact},
\(p^*=\min\{c\mid \Phi(c)=0\}\). Thus the algorithm treats the computation of
\(p^*\) as a scalar root-finding problem for \(\Phi\). To motivate its
Newton-type step,
  suppose that, at a point \(f^*<c<p^*\), an exact global affine minorant of \(\Phi\)
is available:
\[
\Phi(c)-z_c(s-c)\le \Phi(s)
\quad\text{for all }s>f^*,\qquad z_c>0.
\]
The zero of this minorant gives the ideal level update
\(c^+=c+\Phi(c)/z_c\). 
\(\SNtVFA{}\) turns the ideal update above into a safeguarded inexact search.
The initialization first obtains a lower endpoint \(l_0\) through \(\LBInit\)
and a lower-level comparison point \(\bm x_g\) by applying
Lemma~\ref{lem:front-apgz-certificate} with \(\delta_g=\epsilon/8\):
\begin{equation}
\label{eq:newton-reference-xg}
0\le g(\bm x_g)-g^*\le\epsilon/8.
\end{equation}
With \(u_0=f(\bm x_g)\), the algorithm returns \(\bm x_g\) immediately if
\(u_0-l_0\le\epsilon\). Otherwise it starts from \(c_0=l_0\) and tests each
iterate \(c_t\) by calls to \(\VTest\) with decreasing tolerances.
Given a call
\((\bm x_t^j,z_t^j)=\VTest(c_t,\delta_t^j,\bm x_{\rm warm})\),
with \(D_g\) fixed as in~\eqref{eq:Dg-Omegag}, set
\begin{equation}
\mathsf{gap}_t^j:=g(\bm x_t^j)-g(\bm x_g),
\qquad
\mathsf{err}_t^j:=(D_g^2+2)\delta_t^j.
\label{eq:newton-gap-def}
\end{equation}
There are then three possible conclusive cases: a comparison-gap return, a
zero-multiplier return, or a positive-multiplier step. The last case performs
the update 
\[
c_{t+1}=c_t+\frac{\mathsf{gap}_t^j-\mathsf{err}_t^j}{z_t^j}.
\]
  If none of these cases occurs, the same trial level is retested with
  the oracle tolerance halved.

\begin{algorithm}[ht!]
\caption[Safeguarded Newton-type value-function algorithm (\SNtVFA{})]{Safeguarded Newton-type value-function algorithm:
\protect\\
\(\widehat{\bm x} = \SNtVFA(\epsilon,D_g,B_f,
\delta_{\rm init},\bm x_{\rm init}^f,\bm x_{\rm init}^g)\)}
\label{alg:dynamic-newton}
\DontPrintSemicolon
\KwIn{\(D_g\) fixed as in~\eqref{eq:Dg-Omegag},
\(B_f \ge1\) a Lipschitz constant for \(f\) on \({\rm dom}(g_2)\),
\(0<\epsilon\le\min\{1,36\delta_{\rm init}\}\),
\(0<\delta_{\rm init}\le{1/B_f}\),
\(\bm x_{\rm init}^f,\bm x_{\rm init}^g\in{\rm dom}(g_2)\).}
\KwOut{\(\widehat{\bm x}\) with
\(f(\widehat{\bm x})-p^*\le\epsilon\) and
\(g(\widehat{\bm x})-g^*\le\epsilon\).}
Invoke
\((l_0,\Delta,\bm x_f,M_\Delta)\gets
\LBInit(\delta_{\rm init},D_g,B_f,\bm x_{\rm init}^f,\bm x_{\rm init}^g)\).\label{alg:newton-lbinit}\;

Set \(C_D\) as in~\eqref{eq:CD-def} and \(C_\Delta\) as
in~\eqref{eq:CDDelta}.\;
Invoke
\(\bm x_g\gets
\APGz
(g_1,g_2,\epsilon/8,\bm x_{\rm init}^g,D_g)\) \label{alg:newton-xg} and set \(u_0\gets f(\bm x_g)\).\;
\If{\(u_0-l_0\le\epsilon\)\label{alg:newton-early-if}}{
\Return \(\widehat{\bm x}=\bm x_g\).\label{alg:newton-early-return}
\Comment*[r]{early objective-gap return}
}
Set \(c_0\gets l_0\), \(\delta_0^0\gets\delta_{\rm init}\), and \(\bm x_{\rm warm}\gets\bm x_g\).\label{alg:newton-init}\;
\For{\(t=0,1,\ldots\)}{
\For{\(j=0,1,\ldots\)}{
	Invoke
\((\bm x_t^j,z_t^j)\gets
\VTest(c_t,\delta_t^j,\bm x_{\rm warm})\).\label{alg:newton-dual}\;
Set \(\mathsf{gap}_t^j\gets g(\bm x_t^j)-g(\bm x_g)\) and
\(\mathsf{err}_t^j\gets(D_g^2+2)\delta_t^j\).\label{alg:newton-gap-values}\;
\uIf{\(\mathsf{gap}_t^j\le 7\epsilon/8\) and \(\delta_t^j\le\epsilon\)\label{alg:newton-gap-if}}{
\Return \(\widehat{\bm x}=\bm x_t^j\).\label{alg:newton-gap-return}
\Comment*[r]{comparison-gap return}
}
\uElseIf{\(z_t^j=0\) and \(\mathsf{err}_t^j\le\epsilon\)\label{alg:newton-zero-if}}{
\Return \(\widehat{\bm x}=\bm x_t^j\).\label{alg:newton-zero-return}
\Comment*[r]{zero-multiplier return}
	}
	\uElseIf{\(z_t^j>0\) and \(\mathsf{err}_t^j\le \mathsf{gap}_t^j/8\)\label{alg:newton-safe-if}}{
	Set \(c_{t+1}\gets c_t+(\mathsf{gap}_t^j-\mathsf{err}_t^j)/z_t^j\).\label{alg:newton-safe-step}
	\Comment*[r]{positive-multiplier step}
	Set \(\bm x_{\rm warm}\gets\bm x_t^j\) and
	\(\delta_{t+1}^0\gets\delta_t^j\).\label{alg:newton-accept-step}\;
	\Break\;
}
\Else{
Set \(\bm x_{\rm warm}\gets\bm x_t^j\) and \(\delta_t^{j+1}\gets\delta_t^j/2\).\label{alg:newton-refine} \Comment*[r]{ inconclusive comparison step}
}
}
}
\end{algorithm}

The next theorem gives the correctness and operation bound for
\SNtVFA{}. Its proof is given in
Subsection~\ref{subsec:newton-phase-analysis}.
\begin{theorem}[\SNtVFA{} guarantee]
\label{thm:newton-main}
Suppose that Assumption~\ref{ass:basic} holds. Let \(D_g\) be fixed as
in~\eqref{eq:Dg-Omegag}, and let \(B_f \ge1\) be a Lipschitz constant for \(f\) on
\({\rm dom}(g_2)\). Then \SNtVFA{} returns an
\((\epsilon,\epsilon)\)-optimal solution of Problem~\eqref{p:primal}, and, with
\(Q_{\rm init}\) and \(T_{\rm pre}\) defined in \eqref{eq:m0-def} and
Lemma~\ref{lem:front-lbinit}, respectively, has total operation complexity
\[
T_{\rm SNtVFA}
=
T_{\rm pre}
+
\cO\!\left(
D_g\sqrt{\frac{L_{g_1}}{\epsilon}}
+
\sqrt{\frac{C_D(D_g^2+2)(L_{g_1}+L_{f_1}\Omega_g Q_{\rm init})}{\epsilon}}
{ \left(1+\log\frac1\epsilon\right)^3}
\right).
\]

\end{theorem}

This bound has the same interpretation as the discussion following
Theorem~\ref{thm:correctness-bisection}: up to logarithmic factors,
\(T_{\rm SNtVFA}=\tO(\epsilon^{-1/2})\).
\subsection{\texorpdfstring{Newton-Phase Analysis}{Newton-Phase Analysis}}
\label{subsec:newton-phase-analysis}
We collect the standing hypotheses for the Newton-phase analysis in the
following setting.
\begin{setting}[Newton-phase setting]
\label{set:newton-phase}
Assumption~\ref{ass:basic} holds.
The constants \(D_g\) and \(\Omega_g\) are fixed as
in~\eqref{eq:Dg-Omegag}. The constant \(B_f \ge1\) denotes a fixed Lipschitz
constant for \(f\) on \({\rm dom}(g_2)\), whose existence is guaranteed by
Lemma~\ref{lem:lip-bound-f}.
Line~\ref{alg:newton-lbinit} calls \(\LBInit\) to construct a reference certificate
\((l_0,\Delta,\bm x_f,M_\Delta)\) satisfying Condition~\ref{cond:reference},
and Line~\ref{alg:newton-xg} produces a comparison point
\(\bm x_g\) satisfying
\eqref{eq:newton-reference-xg}. For each call to \(\VTest\) in
Line~\ref{alg:newton-dual}, write
\((\bm x_t^j,z_t^j)=\VTest(c_t,\delta_t^j,\bm x_{\rm warm})\).
Here \(0<\delta_t^j\le\delta_{\rm init} {\le 1/B_f}\), and the warm-start point satisfies
\(\bm x_{\rm warm}\in{\rm dom}(g_2)\).
For any executed outer iteration \(t\) and refinement index \(j\) with
\(c_t\ge l_0\), the \(\VTest\) guarantee gives
\begin{equation}
\label{eq:newton-oracle-output}
\bm x_t^j\in {\rm dom}(g_2),\qquad z_t^j\in[0,2M_\Delta),
\end{equation}
and, with \(\mathsf{err}_t^j\) defined in~\eqref{eq:newton-gap-def},
\begin{equation}
\label{eq:newton-oracle-estimates}
f(\bm x_t^j)-c_t\le\delta_t^j,\qquad
{ \lvert g(\bm x_t^j)-\bar g(c_t)\rvert}\le\delta_t^j,\qquad
d_{c_t}(z_t^j)>g(\bm x_t^j)-\mathsf{err}_t^j.
\end{equation}

\end{setting}

The multiplier \(z_t^j\) from \(\VTest\) is used through a global affine
minorant, rather than as a derivative of a smooth value function. The underlying
inequality holds for any \(c\ge l_0\), any \(\tau\in(f^*,+\infty)\), and any
\(z\ge0\):
\begin{equation}
\label{eq:dual-minorant}
d_c(z)-z(\tau-c)
=
\inf_{\bm x\in{\rm dom}(g_2)}\{g(\bm x)+z(f(\bm x)-\tau)\}
\le
\bar g(\tau),
\end{equation}
where the inequality follows by restricting the infimum to
\(\{\bm x\in{\rm dom}(g_2)\mid f(\bm x)\le \tau\}\).

The estimates for the \(\VTest\) call in
\eqref{eq:newton-oracle-estimates}, together with the global minorant
\eqref{eq:dual-minorant}, yield the branch guarantees used in
Algorithm~\ref{alg:dynamic-newton}.
\begin{lemma}[Branch guarantees]
\label{lem:newton-certificates}
In Setting~\ref{set:newton-phase}, consider an outer iteration
\(t\) satisfying \(l_0\le c_t<p^*\), and let \((\bm x_t^j,z_t^j)\) denote the
\(j\)-th oracle output. Then the following assertions hold.

\noindent\textup{(i)} If the condition in Line~\ref{alg:newton-gap-if} or
Line~\ref{alg:newton-zero-if} is satisfied, then \(\bm x_t^j\) is an
\((\epsilon,\epsilon)\)-optimal solution.

\noindent\textup{(ii)} If the condition in Line~\ref{alg:newton-safe-if} is
satisfied, then the update in Line~\ref{alg:newton-safe-step} satisfies
\(c_t<c_{t+1}<p^*\).
\end{lemma}

\begin{proof}
For (i), both return branches give \(f(\bm x_t^j)-p^*\le\epsilon\).
Indeed, \(c_t<p^*\) and the feasibility estimate in
\eqref{eq:newton-oracle-estimates} imply
\(f(\bm x_t^j)-p^*\le f(\bm x_t^j)-c_t\le\delta_t^j\). In the
comparison-gap branch, Line~\ref{alg:newton-gap-if} gives
\(\delta_t^j\le\epsilon\). In the zero-multiplier branch,
Line~\ref{alg:newton-zero-if} gives \(\mathsf{err}_t^j\le\epsilon\), and
\(\delta_t^j\le\mathsf{err}_t^j\) because
\(\mathsf{err}_t^j=(D_g^2+2)\delta_t^j\) and \(D_g^2+2\ge1\).

It remains to bound \(g(\bm x_t^j)-g^*\). If
Line~\ref{alg:newton-gap-if} is satisfied, then
\eqref{eq:newton-reference-xg} gives
\[
g(\bm x_t^j)-g^*
=\mathsf{gap}_t^j+g(\bm x_g)-g^*
\le 7\epsilon/8+\epsilon/8=\epsilon.
\]
If Line~\ref{alg:newton-zero-if} is satisfied, then \(z_t^j=0\) and
\(\mathsf{err}_t^j\le\epsilon\). Since
\(d_{c_t}(0)=\inf_{\bm x\in{\rm dom}(g_2)}g(\bm x)=g^*\), the dual lower
bound in~\eqref{eq:newton-oracle-estimates} gives
\[
g^*=d_{c_t}(0)=d_{c_t}(z_t^j)>
g(\bm x_t^j)-\mathsf{err}_t^j,
\]
and hence \(g(\bm x_t^j)-g^*<\mathsf{err}_t^j\le\epsilon\). Thus either
return branch produces an \((\epsilon,\epsilon)\)-optimal solution.

For (ii), the positive-multiplier branch in Line~\ref{alg:newton-safe-if} gives \(z_t^j>0\) and
\(\mathsf{err}_t^j\le\mathsf{gap}_t^j/8\). Hence
\(\mathsf{gap}_t^j-\mathsf{err}_t^j>0\), and
Line~\ref{alg:newton-safe-step} gives \(c_{t+1}>c_t\). Since
\(c_t\ge l_0>f^*\), we may use \eqref{eq:dual-minorant} with
\((c,z,\tau)=(c_t,z_t^j,c_{t+1})\). Together with
\eqref{eq:newton-oracle-estimates}, this gives
\[
\bar g(c_{t+1})
\ge d_{c_t}(z_t^j)-z_t^j(c_{t+1}-c_t)
> g(\bm x_t^j)-\mathsf{err}_t^j
-\bigl(\mathsf{gap}_t^j-\mathsf{err}_t^j\bigr)
\overset{\eqref{eq:newton-gap-def}}{=}
g(\bm x_g)\overset{\eqref{eq:newton-reference-xg}}{\ge} g^*.
\]
Therefore \(\bar g(c_{t+1})>g^*\). If
\(c_{t+1}\ge p^*\), then Fact~\ref{fact} would imply
\(\bar g(c_{t+1})=g^*\), a contradiction. Hence
\(c_t<c_{t+1}<p^*\).
\end{proof}

Lemma~\ref{lem:newton-certificates} is conditional: once one of the three branch
tests is met, the algorithm either returns an \((\epsilon,\epsilon)\)-optimal
point or accepts a level update that remains below \(p^*\). 
The next lemma first shows that each fixed-level refinement loop exits after
finitely many oracle calls. It then records a uniform lower bound on the
  tolerances used in the \(\VTest\) calls during the Newton phase
and counts all refinements globally.
\begin{lemma}[Inner-refinement guarantee]
\label{lem:newton-inner-finite}
In Setting~\ref{set:newton-phase}, suppose that
\(0<\epsilon\le 36\delta_{\rm init}\) and that every generated level before termination lies in \([l_0,p^*)\). Fix such a level \(c_t\) with
\(l_0\le c_t<p^*\). Then the refinement loop at \(c_t\) exits after finitely
  many \(\VTest\) calls. Moreover, throughout the Newton phase before the algorithm
terminates, every \(\VTest\) call satisfies
\begin{equation}
\label{eq:newton-tolerance-floor}
\delta_t^j\ge
\delta_{\min}^{\rm Nt}:=
\frac{7\epsilon}{128(D_g^2+2)}.
\end{equation}
Let \(J_{\rm ref}^{\rm Nt}\) be the total number of executions of
Line~\ref{alg:newton-refine}. Then
\begin{equation}
\label{eq:newton-global-refinement-bound}
J_{\rm ref}^{\rm Nt}
\le
\left\lceil
\log_2\frac{64(D_g^2+2)\delta_{\rm init}}{7\epsilon}
\right\rceil_+.
\end{equation}
\end{lemma}

\begin{proof}
Set \(J:=J_t^{\rm Nt}:=
\left\lceil\log_2(64(D_g^2+2)\delta_t^0/(7\epsilon))\right\rceil_+\).
By the halving rule in Line~\ref{alg:newton-refine} and the definition of \(J\),
after \(J\) refinements, we have 
\(\mathsf{err}_t^J=(D_g^2+2)\delta_t^J\le 7\epsilon/64\) and
\(\delta_t^J\le \epsilon\).
By \eqref{eq:newton-oracle-output},
\(z_t^J\in[0,2M_\Delta)\). At refinement index \(J\), exactly one of the
following alternatives holds:
\[
\textup{(a)}\ \mathsf{gap}_t^J\le7\epsilon/8,\quad
\textup{(b)}\ \mathsf{gap}_t^J>7\epsilon/8,\ z_t^J=0,\quad
\textup{(c)}\ \mathsf{gap}_t^J>7\epsilon/8,\ z_t^J>0.
\]
Alternative (a) satisfies Line~\ref{alg:newton-gap-if} because
\(\delta_t^J\le\epsilon\). Alternative (b) satisfies
Line~\ref{alg:newton-zero-if} because
\(\mathsf{err}_t^J\le7\epsilon/64<\epsilon\). Under
alternative (c),
\(\mathsf{err}_t^J\le7\epsilon/64<\mathsf{gap}_t^J/8\), so Line~\ref{alg:newton-safe-if}
applies. Hence the refinement loop at \(c_t\) exits no later
than the call with index \(J\), and therefore after finitely many oracle calls.

At the first Newton-phase call, Line~\ref{alg:newton-init} gives
\(\delta_0^0=\delta_{\rm init}\). Since \(D_g^2+2\ge2\) and
\(\epsilon\le36\delta_{\rm init}\),
\(\delta_0^0=\delta_{\rm init}\ge \epsilon/36
>7\epsilon/[128(D_g^2+2)]\).
An accepted step copies its accepted tolerance to the first call at the next
level by Line~\ref{alg:newton-accept-step}, so the lower bound is inherited
across accepted steps. If Line~\ref{alg:newton-refine} is executed after a call
with tolerance \(\delta\), then the fixed-level case analysis above could not have already satisfied one of the branch tests. Hence
\((D_g^2+2)\delta>7\epsilon/64\).
The next tolerance is \(\delta/2\), and therefore it is larger than
\(7\epsilon/[128(D_g^2+2)]\). This proves~\eqref{eq:newton-tolerance-floor}.

Lines~\ref{alg:newton-accept-step} and~\ref{alg:newton-refine}
ensure that the tested Newton-phase tolerances are globally nonincreasing: after an accepted step, its tolerance becomes the initial tolerance for the next outer iteration, and every refinement halves a tolerance \(\delta\) satisfying \((D_g^2+2)\delta>7\epsilon/64\). Starting from
\(\delta_{\rm init}\), the number of such halvings is bounded by
\eqref{eq:newton-global-refinement-bound}.
\end{proof}

It remains to count the accepted positive-multiplier updates.
In the sequel, an \emph{accepted step} means an accepted
positive-multiplier update from Line~\ref{alg:newton-safe-step}. We
first record the residual estimates at the current $c_t$ and affine-minorant
guarantee produced by each accepted step.
\begin{lemma}[Accepted-step guarantee]
\label{lem:newton-affine-minorant}
In Setting~\ref{set:newton-phase}, let \(c_t\in[l_0,p^*)\) be a
level reached by the main loop, and suppose that the oracle call
\((\bm x_t^j,z_t^j)=\VTest(c_t,\delta_t^j,\bm x_{\rm warm})\) at
Line~\ref{alg:newton-dual} satisfies the positive-multiplier condition in
Line~\ref{alg:newton-safe-if}.
Let
\(\kappa=71/49\) and define
\(q:=\mathsf{gap}_t^j-\mathsf{err}_t^j\) and
\(U:=\mathsf{gap}_t^j+\delta_t^j+\epsilon/8\).
Then the following holds:
\begin{equation}
\label{eq:newton-accepted-local}
0<q\le \Phi(c_t)\le U,
\qquad
1\le \frac{U}{q}\le\kappa<2,
\qquad
z_t^j>0,
\qquad
c_{t+1}=c_t+\frac{q}{z_t^j}.
\end{equation}
Moreover, the following affine-minorant estimate holds:
\begin{equation}
\label{eq:newton-affine-minorant-estimate}
q-z_t^j(\tau-c_t)\le \Phi(\tau)
\quad \forall \tau\in(f^*,+\infty).
\end{equation}
\end{lemma}

\begin{proof}
Omit the outer and inner indices and write
\(\bm x=\bm x_t^j\), \(G=\mathsf{gap}_t^j\),
\(E=\mathsf{err}_t^j\), \(\delta=\delta_t^j\), \(z=z_t^j\),
\(q=G-E\), and \(U=G+\delta+\epsilon/8\).

\emph{Residual sandwich.}
Line~\ref{alg:newton-safe-if} gives \(E\le G/8\), and hence \(q>0\).
It also gives \(z>0\). Moreover \(E=(D_g^2+2)\delta\ge\delta\). The value estimate in
\eqref{eq:newton-oracle-estimates}, together with
\eqref{eq:newton-reference-xg}, yields
\[
\Phi(c_t)
=\bar g(c_t)-g^*
\ge g(\bm x)-\delta-g^*
= G+g(\bm x_g)-g^*-\delta
\ge G-\delta
\ge q,
\qquad
\Phi(c_t)\le G+\delta+\epsilon/8=U.
\]

\emph{Relative accuracy.}
If \(G\le7\epsilon/8\), then \(E\le G/8\) and \(E\ge\delta\) would give
\(\delta\le\epsilon\), so Line~\ref{alg:newton-gap-if} would have been satisfied.
Thus \(G>7\epsilon/8\), and hence \(\epsilon/8<G/7\). Together with
\(E\le G/8\), this gives
\[
\frac{U}{q}
=
\frac{G+\delta+\epsilon/8}{G-E}
\le
\frac{G+E+G/7}{G-E}
\le
\frac{G+G/8+G/7}{G-G/8}
=
\frac{71}{49}
=
\kappa.
\]
The lower bound \(U/q\ge1\) follows from
\(0<q\le \Phi(c_t)\le U\).

\emph{Affine minorant.}
The dual lower bound~\eqref{eq:front-dual-lower} and
\(g(\bm x_g)\ge g^*\) give
\(q+g^*=g(\bm x)-g(\bm x_g)-E+g^*
\le g(\bm x)-E<d_{c_t}(z)\).
Combining this estimate with \eqref{eq:dual-minorant} gives, for every
\(\tau\in(f^*,+\infty)\),
\(q+g^*-z(\tau-c_t)<d_{c_t}(z)-z(\tau-c_t)\le \bar g(\tau)\).
This proves \eqref{eq:newton-affine-minorant-estimate}. Since \(z>0\) and
\(q=G-E\), Line~\ref{alg:newton-safe-step} gives \(c_{t+1}=c_t+q/z\).
Together with the preceding estimates, this proves
\eqref{eq:newton-accepted-local}.
\end{proof}

The next lemma combines the complexity result of
\cite[Theorem~2.3]{aravkin2019level} with the comparison-gap test to bound the
number of accepted steps.
\begin{lemma}[Accepted-step bound]
\label{lem:newton-accepted-bound}
In Setting~\ref{set:newton-phase}, assume that, before
termination, each generated level in the Newton-type phase lies in
\([l_0,p^*)\). 
Let \(N_{\rm acc}^{\rm Nt}(\epsilon)\) denote the total number of
accepted steps, namely the accepted positive-multiplier updates
\(c_{t+1}=c_t+(\mathsf{gap}_t^j-\mathsf{err}_t^j)/z_t^j\)
from Line~\ref{alg:newton-safe-step}, in
Algorithm~\ref{alg:dynamic-newton}, and set
\(\kappa:=71/49\). Then
\begin{equation}
\label{eq:newton-accepted-log}
N_{\rm acc}^{\rm Nt}(\epsilon)
\le
\left\lceil
\max\left\{
1+
\log_{2/\kappa}
\frac{
8\max\left\{
2M_\Delta B_fD_g,\,
\Omega_g
\right\}}
{3\epsilon},
2
\right\}
\right\rceil.
\end{equation}
\end{lemma}

\begin{proof}
If no accepted step occurs, then \(N_{\rm acc}^{\rm Nt}(\epsilon)=0\), and the
claim is immediate. Hence assume that at least one accepted step occurs and
index the accepted steps by \(k=0,1,\ldots\). For the \(k\)-th accepted step, write
\((c_k^{\rm acc},\delta_k^{\rm acc},\bm x_k^{\rm acc},z_k^{\rm acc})\)
for the level, tolerance, primal point, and multiplier used at
Line~\ref{alg:newton-safe-step}, and set
\[
q_k:=\mathsf{gap}_k^{\rm acc}-(D_g^2+2)\delta_k^{\rm acc},
\qquad
U_k:=\mathsf{gap}_k^{\rm acc}+\delta_k^{\rm acc}+\epsilon/8,
\qquad
c_k^+:=c_k^{\rm acc}+q_k/z_k^{\rm acc},
\]
where \(\mathsf{gap}_k^{\rm acc}:=g(\bm x_k^{\rm acc})-g(\bm x_g)\).
Before the first accepted step, refinements do not change the level and any
non-accepted inner-loop exit returns from the algorithm. Hence
\(c_0^{\rm acc}=l_0\). By hypothesis, each accepted level, being generated before termination, lies in
\([l_0,p^*)\).

For the first accepted step, set
\(\Theta:=\max\{z_0^{\rm acc}(p^*-c_0^{\rm acc}),q_0\}\).
We first bound the two terms in \(\Theta\). Since
\(c_0^{\rm acc}=l_0\), for any { SBO solution \(\bm x^*\)},
\[
p^*-c_0^{\rm acc}
=p^*-l_0
=f(\bm x^*)-f(\bm x_f)-\Delta
\le { \lvert f(\bm x^*)-f(\bm x_f)\rvert}
\le B_f\|\bm x^*-\bm x_f\|
\le B_fD_g.
\]
The first relation in~\eqref{eq:newton-accepted-local}, applied to the first
accepted step, gives \(q_0\le\Phi(c_0^{\rm acc})=\Phi(l_0)\). Moreover,
\(l_0=f(\bm x_f)+\Delta\) with \(\Delta>0\), so \(\bm x_f\) is feasible
for \({\rm VF}_{l_0}\) and
\(\Phi(l_0)=\bar g(l_0)-g^*\le g(\bm x_f)-g^*\). Combining this estimate
with \(\bm x_f\in{\rm dom}(g_2)\), \eqref{eq:Dg-Omegag}, and
Setting~\ref{set:newton-phase} gives
\(q_0\le \Phi(l_0)\le g(\bm x_f)-g^*\le\Omega_g\) and
\(z_0^{\rm acc}<2M_\Delta\). Consequently,
\(\Theta\le
\bar\Theta:=\max\{2M_\Delta B_fD_g,\,\Omega_g\}\).

For any finite prefix of accepted steps, Fact~\ref{fact} and
\eqref{eq:newton-accepted-local}--\eqref{eq:newton-affine-minorant-estimate}
give the hypotheses of \cite[Theorem~2.3]{aravkin2019level}: \(\Phi\) is
convex and nonincreasing, positive on \((f^*,p^*)\), has smallest root
\(p^*\), and each accepted index satisfies
\[
0<q_k\le \Phi(c_k^{\rm acc})\le U_k,
\qquad
q_k-z_k^{\rm acc}(\tau-c_k^{\rm acc})\le \Phi(\tau)
\quad \forall \tau\in(f^*,+\infty),
\qquad
c_k^+=c_k^{\rm acc}+q_k/z_k^{\rm acc}.
\]
If the next accepted step exists, Line~\ref{alg:newton-safe-step}
sets the next generated level to \(c_k^+\), and refinements do not change this
level; hence \(c_{k+1}^{\rm acc}=c_k^+\).
Moreover, \eqref{eq:newton-accepted-local} in Lemma~\ref{lem:newton-affine-minorant} gives \(U_k/q_k\le\kappa<2\).
Applying this theorem with
\[
K_\epsilon:=
\left\lceil
\max\left\{
1+\log_{2/\kappa}\frac{8\Theta}{3\epsilon},
2
\right\}
\right\rceil
\]
shows that, if the algorithm accepted more than \(K_\epsilon\)
steps, then some accepted index \(k\le K_\epsilon\) would satisfy
\(U_k\le3\epsilon/4\).

For that accepted step,
\(U_k:=\mathsf{gap}_k^{\rm acc}+\delta_k^{\rm acc}+\epsilon/8
\le 3\epsilon/4\).
Line~\ref{alg:newton-safe-if} and \eqref{eq:newton-gap-def} give
\((D_g^2+2)\delta_k^{\rm acc}\le
\mathsf{gap}_k^{\rm acc}/8\). Since
Setting~\ref{set:newton-phase} gives \(\delta_k^{\rm acc}>0\), we have
\(\mathsf{gap}_k^{\rm acc}>0\). The inequality \(U_k\le 3\epsilon/4\) then gives
\(\mathsf{gap}_k^{\rm acc}+\delta_k^{\rm acc}\le5\epsilon/8\), and hence
\(\mathsf{gap}_k^{\rm acc}<5\epsilon/8\) and
\(\delta_k^{\rm acc}<5\epsilon/8\).
Thus the comparison-gap branch in Line~\ref{alg:newton-gap-if} would have been
taken before the positive-multiplier branch in
Line~\ref{alg:newton-safe-if}, which contradicts the fact that this step was
accepted. Therefore \(N_{\rm acc}^{\rm Nt}(\epsilon)\le K_\epsilon\). Since
\(\Theta\le\bar\Theta\), substituting \(\bar\Theta\) into \(K_\epsilon\) gives
\eqref{eq:newton-accepted-log}.
\end{proof}

We are now ready to prove Theorem~\ref{thm:newton-main}. The proof verifies
termination and accuracy first, and then derives the operation bound from the
accepted-step estimate and the value-function oracle cost.
\begin{proof}[Proof of Theorem~\ref{thm:newton-main}]
\emph{Initialization and early exit.}
Lemma~\ref{lem:front-lbinit} gives the reference certificate in
Line~\ref{alg:newton-lbinit}. It satisfies Condition~\ref{cond:reference}, has
preprocessing cost \(T_{\rm pre}\), and obeys
\(M_\Delta\le\Omega_gQ_{\rm init}\). Lemma~\ref{lem:front-apgz-certificate},
applied with \(\delta_g=\epsilon/8\), gives \(\bm x_g\in{\rm dom}(g_2)\) and
\(0\le g(\bm x_g)-g^*\le\epsilon/8\).

If the initial scalar-width test in
Line~\ref{alg:newton-early-if} is satisfied, then \(u_0-l_0\le\epsilon\).
Since \(l_0<p^*\) and \(u_0=f(\bm x_g)\),
\(f(\bm x_g)-p^*=u_0-p^*<u_0-l_0\le\epsilon\) and
\(g(\bm x_g)-g^*\le\epsilon/8\).
Thus the returned point is \((\epsilon,\epsilon)\)-optimal.

\emph{Loop invariant and finite termination.}
Assume now that the Newton-type phase is entered. Before termination, every
generated level satisfies \(l_0\le c_t<p^*\), every tolerance is positive and no
larger than \(\delta_{\rm init}\), and every warm start lies in
\({\rm dom}(g_2)\). These properties hold at the first \(\VTest\) evaluation. Refinement only
halves the tolerance and replaces the warm start by the oracle output in
\({\rm dom}(g_2)\). By
Lemma~\ref{lem:newton-certificates}(ii), an accepted step satisfies
\(c_t<c_{t+1}<p^*\), so \(l_0\le c_{t+1}<p^*\). Line~\ref{alg:newton-accept-step}
sets \(\bm x_{\rm warm}\gets\bm x_t^j\), with
\(\bm x_t^j\in{\rm dom}(g_2)\). Hence
Setting~\ref{set:newton-phase} applies to every call before termination.

By Lemma~\ref{lem:newton-inner-finite}, the inner loop at each generated level
exits after finitely many \(\VTest\) calls. If the algorithm did not terminate,
then every inner-loop exit would have to be an accepted step, since
the other two exits return. This would produce infinitely many accepted steps,
contradicting Lemma~\ref{lem:newton-accepted-bound}. Therefore \(\SNtVFA{}\)
terminates after finitely many iterations. The final exit is either the
comparison-gap branch in Lines~\ref{alg:newton-gap-if}--\ref{alg:newton-gap-return}
or the zero-multiplier branch in
Lines~\ref{alg:newton-zero-if}--\ref{alg:newton-zero-return}, so
Lemma~\ref{lem:newton-certificates}(i) gives an
\((\epsilon,\epsilon)\)-optimal solution.

\emph{Per-call cost and total complexity.}
By Lemma~\ref{lem:front-apgz-certificate}, the comparison APG call in
Line~\ref{alg:newton-xg} costs
\(\cO(D_g\sqrt{L_{g_1}/\epsilon})\).
The Newton phase uses at most
\(N_{\rm acc}^{\rm Nt}(\epsilon)+J_{\rm ref}^{\rm Nt}+1\) \(\VTest\) evaluations,
which is \({ \cO(\log(1/\epsilon))}\) by
\eqref{eq:newton-accepted-log} and \eqref{eq:newton-global-refinement-bound}
for fixed problem data and fixed \(\delta_{\rm init}\).
By~\eqref{eq:newton-tolerance-floor}, every Newton-phase \(\VTest\) call uses
a tolerance at least \(7\epsilon/[128(D_g^2+2)]\). Applying the cost estimate
in Lemma~\ref{lem:front-dual-oracle} with this lower bound gives a uniform
per-call cost. Multiplying it by this call bound, adding the comparison APG cost
and the preprocessing cost \(T_{\rm pre}\), and using
\(M_\Delta\le\Omega_gQ_{\rm init}\), gives the claimed total complexity.
\end{proof}

\section{\texorpdfstring{Value-Function Oracle Construction}{Value-Function Oracle Construction}}
\label{sec:oracle-construction}

This section constructs the value-function oracle used by the
value-function search algorithms and proves Lemma~\ref{lem:front-dual-oracle}.
We first record the APG guarantees and develop the generic
dual-oracle guarantee. We then specialize that guarantee to the perturbed
value-function subproblem.
\subsection{\texorpdfstring{Generic Guarantees and the Dual Oracle}{Generic Guarantees and the Dual Oracle}}
\label{sec:generic-subroutines}

This subsection records the generic first-order guarantees used by the
value-function oracle. For the composite convex problem
\[
\min_{\bm x\in\R^n}~ \varphi(\bm x)\triangleq \varphi_1(\bm x)+\varphi_2(\bm x),
\]
where \(\varphi_1\) is
\(\mu_{\varphi_1}\)-strongly convex \((\mu_{\varphi_1}\ge0)\),
differentiable, and has \(L_{\varphi_1}\)-Lipschitz gradient, and
\(\varphi_2\) is proper, l.s.c., and convex with a computable proximal
mapping, we use two standard APG guarantees: an objective-gap guarantee in the
merely convex case and a stationarity guarantee in the strongly convex case.
We take \(L_{\varphi_1}\) to be a fixed positive upper bound on the
Lipschitz modulus of \(\nabla\varphi_1\), not necessarily the smallest such
bound.

\begin{lemma}[\(\APGz\) objective-gap guarantee {\citep[Theorem~4.4]{beck2009fast}}]
\label{lem:conv-fistalinesearch}
Let \(\delta>0\) and \(\bm x_{\rm init}\in\R^n\). Suppose that
\(D_\varphi\ge{\rm dist}(\bm x_{\rm init},\argmin\varphi)\). A standard
backtracking FISTA call for this composite convex problem with target tolerance \(\delta\),
$
\tilde{\bm x}=\APGz(\varphi_1,\varphi_2,\delta,\bm x_{\rm init},D_\varphi),
$
computes a point \(\tilde{\bm x}\in{\rm dom}(\varphi_2)\) satisfying
\(\varphi(\tilde{\bm x})-\min\varphi\le\delta\) within
\(\cO(D_\varphi\sqrt{L_{\varphi_1}/\delta})\) proximal-gradient operations.
\end{lemma}
\begin{lemma}[\(\APGm\) stationarity guarantee {\citep[Theorem~2.2 and Corollary~2.3]{xu2022first}}]
\label{lem:apgmu-complexity}
Suppose that \({\rm dom}(\varphi_2)\) has diameter at most \(D_\varphi\).
Let \(\mu_{\varphi_1}>0\), \(\delta>0\), and
\(\bm x_{\rm init}\in{\rm dom}(\varphi_2)\). A strongly convex APG call with
target tolerance \(\delta\),
$	
\tilde{\bm x}=\APGm(\varphi_1,\varphi_2,\delta,\bm x_{\rm init}),
$
returns \(\tilde{\bm x}\in{\rm dom}(\varphi_2)\) satisfying
$
{\rm dist}(\bm 0,\partial\varphi(\tilde{\bm x}))\le\delta
$
using
$\cO\!\left(
\sqrt{\frac{L_{\varphi_1}}{\mu_{\varphi_1}}}
\left(1+
\left\lceil\log_2\frac{L_{\varphi_1} D_\varphi}{\delta}\right\rceil_+
\right)
\right)
$
proximal-gradient operations.
\end{lemma}

We next formulate the constrained composite problem solved by the generic dual
oracle.
The generic problem is
\begin{equation}
\label{p:func-cons}
\min\limits_{\bm x\in \R^n} \phi(\bm x)\triangleq \phi_1(\bm x) + \phi_2(\bm x) \qquad
{\rm s.t.} \quad\psi(\bm x)\triangleq \psi_1(\bm x) + \psi_2(\bm x) \le 0.
\end{equation}
For \(z\ge0\), define
\begin{equation}\label{p:dzphi}
\cL_{\phi}(\bm x,z)\triangleq\phi(\bm x)+z\psi(\bm x),
\qquad
d_\phi(z)\triangleq \min_{\bm x\in\R^n}\cL_\phi(\bm x,z).
\end{equation}
A dual solution is any
\(\bar z_\phi\in\argmax_{z\ge0}d_\phi(z)\).
\begin{assumption}
\label{ass:dual-subproblem}
Assume that, for a prescribed bound \(M_z^\phi>0\), the following conditions hold.

\noindent\textup{(i)} \(\phi_1,\psi_1:\R^n\to\R\) are convex and continuously differentiable;
\(\phi_1\) is \(\mu_{\phi_1}\)-strongly convex with
\(\mu_{\phi_1}>0\), and \(\nabla\phi_1\) and \(\nabla\psi_1\) are
\(L_{\phi_1}\)- and \(L_{\psi_1}\)-Lipschitz continuous, respectively.

\noindent\textup{(ii)} \(\phi_2,\psi_2:\R^n\to\R\cup\{+\infty\}\) are proper, l.s.c., and convex,
with \({\rm dom}(\phi_2)\subseteq{\rm dom}(\psi_2)\). The function
\(\psi_2\) is \(l_{\psi_2}\)-Lipschitz continuous on
\({\rm dom}(\phi_2)\).

\noindent\textup{(iii)} The set \({\rm dom}(\phi_2)\) is compact,
with \(D_\phi\ge{\rm diam}({\rm dom}(\phi_2))\) denoting a fixed diameter bound.

\noindent\textup{(iv)} There exists a dual solution
\(\bar z_\phi\in\argmax_{z\ge0}d_\phi(z)\) satisfying $0\le \bar z_\phi\le M_z^\phi$.
\end{assumption}
We use the same operation-counting convention as in
Definition~\ref{def:basic-operation} for the generic pair \((\phi,\psi)\).
Thus component and gradient evaluations, together with the combined proximal
calls for \(\phi_2+z\psi_2\) made by the oracle, count as basic operations.
For \(z\in[0,2M_z^\phi)\), the gradient of the smooth part
\(\phi_1+z\psi_1\) is Lipschitz continuous with constant at most
\(L_{\phi_1}+2M_z^\phi L_{\psi_1}\).
For every \(z\in[0,2M_z^\phi)\),
\(\cL_\phi(\cdot,z)\) is proper, l.s.c., and \(\mu_{\phi_1}\)-strongly convex on
the compact effective domain \({\rm dom}(\phi_2)\). Hence it has a unique
minimizer
\(\bm x_\phi(z)\triangleq \argmin_{\bm x \in \R^n} \cL_{\phi}(\bm x,z)\).

Before presenting the generic dual oracle, we specify the KKT residual bounds returned by the oracle.
\begin{definition}
Given \(\varepsilon_1,\varepsilon_2,\varepsilon_3 \ge 0\), a pair
\((\tilde{\bm x},\tilde z)\) with
\(\tilde{\bm x}\in{\rm dom}(\phi_2)\) and \(\tilde z\ge0\) is an
\((\varepsilon_1,\varepsilon_2,\varepsilon_3)\)-KKT pair of
Problem~\eqref{p:func-cons} if it satisfies the following stationarity,
feasibility, and complementarity bounds:
\begin{equation}
\label{eq:kkt-tol}
{\rm dist}(\bm 0,\partial_{\bm x}\cL_{\phi}(\tilde{\bm x},\tilde z))
\le \varepsilon_1,\qquad
[\psi(\tilde{\bm x})]_+
\le \varepsilon_2,\qquad
\lvert\tilde z\psi(\tilde{\bm x})\rvert
\le \varepsilon_3.
\end{equation}
\end{definition}
Let \(B_\psi{\ge1}\) be any known Lipschitz constant for \(\psi\) on
\({\rm dom}(\phi_2)\). Its existence is recorded in
Lemma~\ref{lem:lip-psi}.

For a dual-oracle tolerance \(\delta\), we use the local shorthand
\begin{equation}
\label{eq:dual-tolerances}
\varepsilon_1:=\mu_{\phi_1}\delta^2,
\qquad
\varepsilon_2:=\delta,
\qquad
\varepsilon_3:=4M_z^\phi\delta.
\end{equation}
When \(\delta\) is fixed by the current routine, we write
\(\varepsilon_i\) for the quantities in~\eqref{eq:dual-tolerances}.

\begin{algorithm}[htp!]
\caption[Dual oracle (\Dual{})]{Dual oracle:
\((\tilde{\bm x},\tilde z)=\Dual(\phi,\psi,\sigma,
M_z^\phi,\bm x_{\rm init};\delta)\)}
\label{alg:dual-algo}
\DontPrintSemicolon
\KwIn{\((\phi,\psi)\) satisfies Assumption~\ref{ass:dual-subproblem},
\(0<\sigma\le M_z^\phi\), \(0<\delta\le1/B_\psi\), and
\(\bm x_{\rm init}\in{\rm dom}(\phi_2)\).}
Invoke \(({\rm flag},\tilde{\bm x},\tilde z,Z,\tilde{\bm x}_a,\tilde{\bm x}_b)
=\IntV(\phi,\psi,\sigma,
M_z^\phi,\bm x_{\rm init};\delta)\).\;
\If{\({\rm flag} = 0\)}{
Invoke \((\tilde{\bm x},\tilde z){\gets}
\Bisec
(\phi, \psi, Z, M_z^\phi, \tilde{\bm x}_a, \tilde{\bm x}_b,
\bm x_{\rm init};\delta)\).\;
}
\Return \((\tilde{\bm x},\tilde z)\).
\end{algorithm}

Algorithm~\ref{alg:dual-algo} is the generic dual oracle for
Problem~\eqref{p:func-cons}, combining the interval search \(\IntV\) with the
multiplier bisection routine \(\Bisec\). 
The multiplier search adapts the interval-search and bisection idea for a
single functional constraint~\citep{xu2022first}. Appendix~\ref{app:dual-oracle-details}
proves the component guarantees for \(\IntV\) and \(\Bisec\) and then composes
them to establish the generic dual-oracle guarantee below. This guarantee
yields the three KKT residual bounds used later in the
value-function oracle.

\begin{proposition}[Generic dual-oracle guarantee]
\label{prop:conv-dual-algo}
Suppose that Assumption~\ref{ass:dual-subproblem} holds for
Problem~\eqref{p:func-cons}. Then the call
\(\Dual(\phi,\psi,\sigma,M_z^\phi,\bm x_{\rm init};\delta)\)
returns an \((\varepsilon_1,\varepsilon_2,\varepsilon_3)\)-KKT pair \((\tilde{\bm x},\tilde z)\) with
\(\tilde z\in[0,2M_z^\phi)\), where
\(\varepsilon_1,\varepsilon_2,\varepsilon_3\) are given by~\eqref{eq:dual-tolerances}.
Moreover, its operation complexity is bounded by
\[
\cO\!\left(
\sqrt{\frac{L_{\phi_1}+M_z^\phi L_{\psi_1}}{\mu_{\phi_1}}}
\left(
1+
\left\lceil
\log_2\frac{(L_{\phi_1}+M_z^\phi L_{\psi_1})D_{\phi}}{\mu_{\phi_1}\delta^2}
\right\rceil_+
\right)
\left[
1+
\left\lceil\log_2\frac{M_z^\phi}{\sigma}\right\rceil_+
+
\left\lceil
\log_2\frac{M_z^\phi}{\mu_{\phi_1}\delta^3}
\right\rceil_+
\right]
\right).
\]
\end{proposition}

For fixed subproblem parameters and fixed \(\sigma\), the displayed bound is
{\(\cO(({ \log(1/\delta)})^2)\) as \(\delta\downarrow0\)}.
Thus, at the oracle level, the dual-based construction removes one logarithmic
factor from the inexact augmented-Lagrangian construction of
\citep[Algorithm~4]{xu2022first} in the smooth-constraint setting, whose
corresponding dependence is cubic.

\subsection{\texorpdfstring{Value-Function Oracle Specialization}{Value-Function Oracle Specialization}}
\label{sec:pre-dual}

For a scalar \(c\), we apply Proposition~\ref{prop:conv-dual-algo} to a
strongly convex perturbation of~\eqref{p:subp-rewrite}. Let
\(\bm x_f\in{\rm dom}(g_2)\) and \(\Delta>0\) satisfy
\(f_c(\bm x_f)\le-\Delta\). Using the Lipschitz constant \(B_f\) fixed in
Lemma~\ref{lem:lip-bound-f}, choose \(0<\eta\le{1/B_f}\) and set
\(G_{\eta}(\bm x):=g(\bm x)+\frac{\eta}{2}\|\bm x-\bm x_f\|^2\).
The resulting perturbed constrained value problem is
\begin{equation}
\label{p:perturbed-str}\tag{$\rm VF_c^\eta$}
\min_{\bm x\in\R^n}\ G_{\eta}(\bm x)
\qquad
{\rm s.t.}\quad
f_c(\bm x)\le 0.
\end{equation}

The following preparatory lemma provides the two estimates used in this
specialization. Part~\textup{(i)} uses the strict-feasibility margin
\(f_c(\bm x_f)\le-\Delta\) to obtain a uniform multiplier bound, and
Part~\textup{(ii)} converts an approximate KKT pair for the perturbed problem
into two-sided value control for \(\bar g(c)\).
\begin{lemma}[Value-function specialization estimates]
\label{lem:bound-mul-per}
Suppose that Assumption~\ref{ass:basic} holds,
\(0<\eta {\le 1/B_f}\), \(\Delta>0\), and
\(\bm x_f\in{\rm dom}(g_2)\) satisfies
\(c\ge f(\bm x_f)+\Delta\). Let
\(\bm x_g^{\rm ref}\in {\rm dom}(g_2)\) satisfy
\(g(\bm x_g^{\rm ref})-g^*\le 1\).
Let \(M_\Delta\) be given by the defining equality in~\eqref{eq:Mz-def}.
Then the following assertions hold.

\noindent\textup{(i)} Problems~\eqref{p:subp-rewrite} and~\eqref{p:perturbed-str}
admit primal-dual solutions. Every optimal multiplier
\(z_c^*\) of Problem~\eqref{p:subp-rewrite} and every optimal multiplier
\(z_\eta^*\) of Problem~\eqref{p:perturbed-str} satisfy
\(\max\{z_c^*,z_\eta^*\}\le M_\Delta\).

\noindent\textup{(ii)} Let \(\tilde{\bm x}\in{\rm dom}(g_2)\), and let
\(\tilde z\ge0\) satisfy
\begin{equation}
\label{eq:kkt-res}
{\rm dist}\bigl(\bm 0,\partial_{\bm x} (G_\eta+\tilde z f_c)(\tilde{\bm x})\bigr)\le \varepsilon_1,\quad
[f_c(\tilde{\bm x})]_+\le \varepsilon_2,\quad
{ \lvert\tilde z f_c(\tilde{\bm x})\rvert}\le \varepsilon_3
\end{equation}
for some \(\varepsilon_1,\varepsilon_2,\varepsilon_3\ge0\). Then,
with \(D_g\) denoting the diameter bound fixed in~\eqref{eq:Dg-Omegag}, we have
\begin{equation}
\label{eq:kkt-value-est}
-M_\Delta \varepsilon_2
\le
g(\tilde{\bm x})-\bar g(c)
\le
D_g \varepsilon_1+\varepsilon_3+\frac{\eta}{2}D_g^2.
\end{equation}
\end{lemma}

\begin{proof}

\textup{(i)} All feasibility and multiplier statements are understood relative to the compact
convex domain \({\rm dom}(g_2)\). The point \(\bm x_f\) satisfies
\(f_c(\bm x_f)\le-\Delta\), so both feasible sets are nonempty. Compactness of
\({\rm dom}(g_2)\), lower semicontinuity of the objectives, and continuity of
\(f\) on \({\rm dom}(g_2)\) give primal solutions. Since \({\rm dom}(g_2)\) is nonempty and convex in \(\mathbb R^n\),
its relative interior \({\rm ri}({\rm dom}(g_2))\) is nonempty. Choose
\(\bm x_{\rm ri}\in{\rm ri}({\rm dom}(g_2))\), and let
$
\bm x_\alpha=(1-\alpha)\bm x_f+\alpha \bm x_{\rm ri}.
$
For all sufficiently small \(\alpha>0\),
\(\bm x_\alpha\in{\rm ri}({\rm dom}(g_2))\) and \(f_c(\bm x_\alpha)<0\).
Hence the relative Slater condition holds. 
By \cite[Theorem~28.2]{Rockafellar1970convex}, both
problems~\eqref{p:subp-rewrite} and~\eqref{p:perturbed-str} admit
primal-dual solutions.

Let \((\bm x_\eta^*,z_\eta^*)\) be any primal-dual solution of
the perturbed problem. 
The KKT conditions give \(z_\eta^*\ge0\),
\(z_\eta^*f_c(\bm x_\eta^*)=0\), and
\(\bm x_\eta^*\in\argmin_{\bm x\in{\rm dom}(g_2)}
\{G_\eta(\bm x)+z_\eta^*f_c(\bm x)\}\). Hence
\[
G_\eta(\bm x_\eta^*)
=G_\eta(\bm x_\eta^*)+z_\eta^*f_c(\bm x_\eta^*)
\le G_\eta(\bm x_f)+z_\eta^*f_c(\bm x_f).
\]
Since \(c\ge f(\bm x_f)+\Delta\), we have
\(f_c(\bm x_f)\le-\Delta\). Together with \(z_\eta^*\ge0\), this gives
\[
G_\eta(\bm x_\eta^*)\le G_\eta(\bm x_f)+z_\eta^*f_c(\bm x_f)
\le G_\eta(\bm x_f)-z_\eta^*\Delta.
\]
The relations \(G_\eta(\bm x_f)=g(\bm x_f)\), \(G_\eta\ge g\), and
\(g(\bm x_g^{\rm ref})-g^*\le1\) give
\[
z_\eta^*\le\frac{G_\eta(\bm x_f)-G_\eta(\bm x_\eta^*)}{\Delta}
\le\frac{g(\bm x_f)-g^*}{\Delta}
\le\frac{g(\bm x_f)-g(\bm x_g^{\rm ref})+1}{\Delta}\le M_\Delta.
\]
The last inequality uses~\eqref{eq:Mz-def}. The same argument with \(\eta=0\)
and an arbitrary primal-dual solution gives the bound for \(z_c^*\).

\textup{(ii)} Let \((\bm x_c^*,z_c^*)\) be a primal-dual solution of
Problem~\eqref{p:subp-rewrite}. By~\eqref{eq:kkt-res}, choose
\(\xi\in\partial_{\bm x}(G_\eta+\tilde z f_c)(\tilde{\bm x})\) with
\(\|\xi\|\le\varepsilon_1\). 
Since \(\tilde{\bm x},\bm x_c^*\in{\rm dom}(g_2)\), the diameter bound gives
\(\|\tilde{\bm x}-\bm x_c^*\|\le D_g\).
The subgradient inequality for \(G_\eta+\tilde z f_c\) at
\(\tilde{\bm x}\), together with feasibility of \(\bm x_c^*\) and
\(\tilde z\ge0\), gives
\[
G_\eta(\tilde{\bm x})-G_\eta(\bm x_c^*)
\le \tilde z\bigl(f_c(\bm x_c^*)-f_c(\tilde{\bm x})\bigr)
   +\langle \xi,\tilde{\bm x}-\bm x_c^*\rangle \le { \lvert\tilde z f_c(\tilde{\bm x})\rvert}+D_g\|\xi\|
\le\varepsilon_3+D_g\varepsilon_1.
\]

The bounds \(G_\eta\ge g\) and
\(G_\eta(\bm x_c^*)\le g(\bm x_c^*)+\eta D_g^2/2\) give the second
inequality in~\eqref{eq:kkt-value-est}.
For the lower bound, Lagrangian optimality and complementarity for
\((\bm x_c^*,z_c^*)\) imply
\[
g(\tilde{\bm x})-\bar g(c)\ge -z_c^*[f_c(\tilde{\bm x})]_+
\ge -M_\Delta\varepsilon_2,
\]
where \eqref{eq:kkt-res} gives \([f_c(\tilde{\bm x})]_+\le\varepsilon_2\), and
Part~\textup{(i)} gives \(z_c^*\le M_\Delta\).
\end{proof}

We now prove the value-function oracle guarantee by applying
Proposition~\ref{prop:conv-dual-algo} to the specialized perturbed
problem~\eqref{p:perturbed-str}.

\begin{proof}[Proof of Lemma~\ref{lem:front-dual-oracle}]
Set \[\delta_{\rm D}:=\delta/C_\Delta, \quad\eta:=\delta/C_D,\quad \text{ and }\quad
\cL_c^\eta(\bm x,z):=G_\eta(\bm x)+z f_c(\bm x).\]
Condition~\ref{cond:reference} gives
\(c\ge l_0=f(\bm x_f)+\Delta\), \(\Delta>0\), and a point
\(\bm x_g^{\rm ref}\in{\rm dom}(g_2)\) with
\(g(\bm x_g^{\rm ref})-g^*\le1\). Hence
Lemma~\ref{lem:bound-mul-per}\textup{(i)},
together with Assumption~\ref{ass:basic} and Lemma~\ref{lem:lip-bound-f},
shows that the decomposition
\[
\phi_1=g_1+\frac{\eta}{2}\|\cdot-\bm x_f\|^2,\quad
\phi_2=g_2,\quad
\psi_1=f_1-c,\quad
\psi_2=f_2
\]
fits Assumption~\ref{ass:dual-subproblem} for
\((\phi,\psi)=(G_\eta,f_c)\), with
\begin{equation}
\label{eq:front-dual-param-match}
\mu_{\phi_1}=\eta,\qquad D_\phi=D_g,\qquad B_\psi=B_f,\qquad
M_z^\phi=M_\Delta.
\end{equation}
By~\eqref{eq:Mz-def}, \(M_\Delta\ge1\). Since the specialized call
\eqref{eq:dynamic-dual-call} passes \(\sigma=1\) to \(\Dual\), the
input condition \(0<\sigma\le M_z^\phi\) is satisfied. Moreover,
\[
0<\delta_{\rm D}\le \delta\le {1/B_f},
\qquad
0<\eta\le {1/B_f}.
\]
Thus the call in~\eqref{eq:dynamic-dual-call} is admissible for
Proposition~\ref{prop:conv-dual-algo}.
The proposition returns
\(\bm x_c\in{\rm dom}(g_2)\) and \(z_c\in[0,2M_\Delta)\), and gives
the following KKT residual bounds:
\begin{equation}
\label{eq:front-dual-kkt-res}
{\rm dist}\bigl(\bm 0,\partial_{\bm x}\cL_c^\eta(\bm x_c,z_c)\bigr)
\le\eta\delta_{\rm D}^2,\quad
[f_c(\bm x_c)]_+\le\delta_{\rm D},\quad
{ \lvert z_c f_c(\bm x_c)\rvert}\le4M_\Delta\delta_{\rm D}=\eta.
\end{equation}
Since \(\delta_{\rm D}=\delta/C_\Delta\) and
\(C_\Delta=4C_DM_\Delta\ge1\), we have
\(\delta_{\rm D}\le\delta\).
The second bound in~\eqref{eq:front-dual-kkt-res} gives
\(f_c(\bm x_c)\le\delta\). 
Applying the estimate~\eqref{eq:kkt-value-est} in
Lemma~\ref{lem:bound-mul-per}\textup{(ii)} to the same KKT pair gives
\[
-M_\Delta\delta_{\rm D}\le g(\bm x_c)-\bar g(c)
\le D_g\eta\delta_{\rm D}^2+\eta+\frac{\eta D_g^2}{2}.
\]
Using \(\delta_{\rm D}\le1\), \(C_D=1+D_g(1+D_g)\), and
\(C_\Delta=4C_DM_\Delta\), the right-hand side is at most \(\delta\) and
the lower bound is at least \(-\delta/(4C_D)\). This proves
\eqref{eq:front-dual-value-cert}.

For the lower bound~\eqref{eq:front-dual-lower}, let
\(d_c^\eta(z):=\inf_{\bm x\in{\rm dom}(g_2)}\cL_c^\eta(\bm x,z)\).
The \(\eta\)-strong convexity of \(\cL_c^\eta(\cdot,z_c)\),
the first bound in~\eqref{eq:front-dual-kkt-res}, and
\(\delta_{\rm D}\le1\) imply
\[
\cL_c^\eta(\bm x_c,z_c)-d_c^\eta(z_c)
\le \frac{1}{2\eta}
{\rm dist}(\bm 0,\partial_{\bm x}\cL_c^\eta(\bm x_c,z_c))^2
\le\frac{\eta}{2}.
\]
Together with the third bound in~\eqref{eq:front-dual-kkt-res}, and using
\(0\le G_\eta(\bm x)-g(\bm x)\le\eta D_g^2/2\) on
\({\rm dom}(g_2)\), this yields the chain
\begin{align*}
d_c(z_c)
&\ge d_c^\eta(z_c)-\frac{\eta D_g^2}{2}
\ge G_\eta(\bm x_c)-\frac{D_g^2+3}{2}\eta
\ge g(\bm x_c)-\frac{D_g^2+3}{2}\eta.
\end{align*}
The relations \(\eta\le\delta\) and
\((D_g^2+3)/2<D_g^2+2\) yield~\eqref{eq:front-dual-lower}.

It remains to insert these parameters into the cost estimate of
Proposition~\ref{prop:conv-dual-algo}, applied here with tolerance
\(\delta_{\rm D}\). In addition to
\eqref{eq:front-dual-param-match}, we use
\[
L_{\phi_1}=L_{g_1}+\eta,\qquad L_{\psi_1}=L_{f_1},\qquad
\sigma=1,\qquad L_{\phi_1}+M_z^\phi L_{\psi_1}
\le 2\bigl(L_{g_1}+M_\Delta L_{f_1}\bigr).
\]
The displayed inequality follows from \(\eta\le1\) and the convention
\(L_{g_1}\ge1\). Since
\(\delta_{\rm D}=\delta/C_\Delta\), for fixed problem and
reference tuple data the two logarithmic factors in
Proposition~\ref{prop:conv-dual-algo} are both
{\(\cO({ \log(1/\delta)})\) as \(\delta\downarrow0\)}. Hence the claimed cost bound
\eqref{eq:front-dual-cost} follows.
\end{proof}

\section{\texorpdfstring{Initialization of the Reference Certificate}{Initialization of the Reference Certificate}}
\label{sec:auxiliary}

This section presents \(\LBInit{}\) (Algorithm~\ref{alg:lb-initialization}) for
constructing the reference certificate used in Lemma~\ref{lem:front-lbinit}. 
It repeatedly performs two APG calls with a \(\Dual\)-based
value-function test at the trial level \(l=f(\bm x_f)+\Delta\), while
halving \(\Delta\) until a strict lower endpoint is verified.
The routine returns
\((l_0,\Delta,\bm x_f,M_\Delta)\) satisfying
\begin{equation}
\label{eq:init-target-tuple}
f^*<l_0<p^*,
\qquad
l_0=f(\bm x_f)+\Delta,
\qquad
\Delta>0,
\qquad
M_\Delta\hbox{ as in }\eqref{eq:Mz-def}.
\end{equation}

\begin{algorithm}[H]
\caption[Lower-endpoint initialization (\LBInit{})]{Lower-endpoint initialization:
\((l,\Delta,\bm x_f,M_\Delta)=\LBInit(\delta_{\rm init},D_g,B_f,
\bm x_{\rm init}^f,\bm x_{\rm init}^g)\)}
\label{alg:lb-initialization}
\DontPrintSemicolon
\KwIn{\(D_g>0\), \(B_f \ge1\),
\(\bm x_{\rm init}^f,\bm x_{\rm init}^g\in{\rm dom}(g_2)\), and
\(0<\delta_{\rm init}\le {1/B_f}\).}
Set \(\Delta\gets2\delta_{\rm init}\),
\(C_D\) as in~\eqref{eq:CD-def},
\(\bm x_f\gets\bm x_{\rm init}^f\),
\(\bm x_g^{\rm ref}\gets\bm x_{\rm init}^g\), and
\(\tilde{\bm x}\gets\bm x_{\rm init}^f\).\; \label{alg:lb-init-params}
\For{\(k=0,1,2,\ldots\)}{
\(\Delta\gets\Delta/2\).\label{alg:lb-update-delta}\Comment*[r]{trial margin refinement}
Invoke
\(\bm x_f\gets
\APGz
(f_1,f_2+{\rm I}_{{\rm dom}(g_2)},\Delta,\bm x_f,D_g)\).\label{alg:lb-apgf}\;
Invoke
\(\bm x_g^{\rm ref}\gets
\APGz
(g_1,g_2,\Delta,\bm x_g^{\rm ref},D_g)\).\label{alg:lb-apgg}
\Comment*[r]{APG comparison points}
Set \(l\gets f(\bm x_f)+\Delta\).\label{alg:lb-build}
\Comment*[r]{trial endpoint}
Set \(M_\Delta\) by~\eqref{eq:Mz-def}, \(C_\Delta\) by~\eqref{eq:CDDelta},
\(\eta\gets\Delta/C_D\), and
\(G_\eta(\bm x)\gets g(\bm x)+\frac{\eta}{2}\|\bm x-\bm x_f\|^2\).\label{alg:lb-compute-const}\;
Invoke \((\tilde{\bm x},\cdot)\gets
\Dual(G_\eta,f-l,1,M_\Delta,\tilde{\bm x};\Delta/C_\Delta)\).
\label{alg:lb-pre-dualcall}
\Comment*[r]{value-function check}
\If{\(g(\tilde{\bm x})-g(\bm x_g^{\rm ref})>\Delta\)}{\label{alg:lb-check}
\Return \((l,\Delta,\bm x_f,M_\Delta)\). \label{alg:lb-return}
	\Comment*[r]{reference tuple \eqref{eq:init-target-tuple}}
}
}
\end{algorithm}

We next record the per-iteration estimates that justify the stopping test in
Line~\ref{alg:lb-check}. Denote by
\(\Delta_k:=\delta_{\rm init}/2^k\) the value of \(\Delta\) after
Line~\ref{alg:lb-update-delta} in iteration \(k\), and let
\(\bm x_f^{(k)}\), \(\bm x_g^{(k)}\), \(l^{(k)}\),
\(M^{(k)}\), and \(\tilde{\bm x}^{(k)}\) be the
corresponding quantities generated in
Lines~\ref{alg:lb-apgf}--\ref{alg:lb-pre-dualcall}. The two APG
calls and Line~\ref{alg:lb-build} give
\begin{align}
&
f(\bm x_f^{(k)})\le f^*+\Delta_k,\qquad
g(\bm x_g^{(k)})\le g^*+\Delta_k,\qquad
\bm x_f^{(k)},\bm x_g^{(k)}\in{\rm dom}(g_2),
\label{eq:init-loop-apg-est}\\
&
\Delta_k>0,\qquad
l^{(k)}=f(\bm x_f^{(k)})+\Delta_k,\qquad
f^*<l^{(k)}\le f^*+2\Delta_k.
\label{eq:init-loop-endpoint-est}
\end{align}

The second relation in~\eqref{eq:init-loop-endpoint-est} gives
\(f(\bm x_f^{(k)})-l^{(k)}=-\Delta_k\), so \(\bm x_f^{(k)}\) is strictly feasible for
\(f(\bm x)-l^{(k)}\le0\) with margin \(\Delta_k\). Also,
\(\Delta_k\le\delta_{\rm init}\le{1/B_f}\), and the \(g\)-estimate
in~\eqref{eq:init-loop-apg-est} gives
\(g(\bm x_g^{(k)})-g^*\le\Delta_k\le1\). Hence, with
\((c,\Delta,\bm x_f,\bm x_g^{\rm ref})
=(l^{(k)},\Delta_k,\bm x_f^{(k)},\bm x_g^{(k)})\),
let \(M^{(k)}\) be the value of \(M_\Delta\) under this substitution, and set
\(\eta^{(k)}:=\Delta_k/C_D\). 
Then Lemma~\ref{lem:bound-mul-per}\textup{(i)} gives, in particular, the multiplier
bound \(z_{\eta^{(k)}}^*\le M^{(k)}\) for Problem~\eqref{p:perturbed-str}. Hence
Assumption~\ref{ass:dual-subproblem}\textup{(iv)} holds for the dual-oracle
call in Line~\ref{alg:lb-pre-dualcall}, with \(M_z^\phi=M^{(k)}\).
Equation~\eqref{eq:Mz-def} gives \(M^{(k)}\ge1\). Since Line~\ref{alg:lb-pre-dualcall}
passes \(\sigma=1\) to \(\Dual\), the input requirement
\(0<\sigma\le M_z^\phi\) is satisfied.

Consequently,
Line~\ref{alg:lb-compute-const} sets \(C_{\Delta}^{(k)}=4C_DM^{(k)}\) and
\(\eta=\eta^{(k)}\). The initialization in
Line~\ref{alg:lb-init-params} gives
\(\tilde{\bm x}=\bm x_{\rm init}^f\in{\rm dom}(g_2)\), so the first call in
Line~\ref{alg:lb-pre-dualcall} has a valid warm start. 
Proposition~\ref{prop:conv-dual-algo} returns an approximate KKT pair
\((\tilde{\bm x}^{(k)},\tilde z^{(k)})\) for Problem~\eqref{p:perturbed-str} with
\(c=l^{(k)}\) and \(\eta=\eta^{(k)}\). 
Applying the estimate~\eqref{eq:kkt-value-est} in
Lemma~\ref{lem:bound-mul-per}\textup{(ii)}
then converts this KKT pair into the value-function estimates
\begin{equation}
\label{eq:init-loop-vtest-est}
{ \lvert g(\tilde{\bm x}^{(k)})-\bar g(l^{(k)})\rvert}\le\Delta_k,
\qquad
f(\tilde{\bm x}^{(k)})-l^{(k)}\le\Delta_k,
\qquad
\tilde{\bm x}^{(k)}\in{\rm dom}(g_2),
\end{equation}
and \(\tilde{\bm x}^{(k)}\in{\rm dom}(g_2)\)
further makes the returned point a valid warm start
for the next iteration.

Thus we claim that the test in Line~\ref{alg:lb-check} is satisfied whenever
\(\Delta_k\le\theta:=\min\{(p^*-f^*)/4,m_0/4\}\). Indeed, if \(\Delta_k\le\theta\), then
\eqref{eq:init-loop-endpoint-est} gives \(l^{(k)}\le(p^*+f^*)/2\). 
By the definition of \(m_0\) in~\eqref{eq:m0-def} and the monotonicity of \(\bar g\), we have
\(\bar g(l^{(k)})\ge g^*+m_0\). 
This yields
\begin{equation}
\label{eq:init-loop-trigger}
g(\tilde{\bm x}^{(k)})-g(\bm x_g^{(k)})
\overset{(\ref{eq:init-loop-apg-est},\ref{eq:init-loop-vtest-est})}{\ge} \bar g(l^{(k)})-\Delta_k-(g^*+\Delta_k)
\ge m_0-2\Delta_k
\ge 2\Delta_k
>\Delta_k.
\end{equation}
By~\eqref{eq:init-loop-trigger}, the test in Line~\ref{alg:lb-check} is
satisfied whenever \(\Delta_k\le\theta\). Since
\(\Delta_k=\delta_{\rm init}/2^k\downarrow0\),
Algorithm~\ref{alg:lb-initialization} terminates. Let \(K_0\) denote its
stopping iteration, and set
\[
l:=l^{(K_0)},\qquad
\Delta:=\Delta_{K_0},\qquad
\bm x_f:=\bm x_f^{(K_0)},\qquad
\bm x_g^{\rm ref}:=\bm x_g^{(K_0)},\qquad
M_\Delta:=M^{(K_0)}.
\]
At this iteration, the test in Line~\ref{alg:lb-check}
and~\eqref{eq:init-loop-vtest-est} imply
\(\bar g(l)>g(\bm x_g^{\rm ref})\ge g^*\). Together with
\(f^*<l\) from~\eqref{eq:init-loop-endpoint-est}, Fact~\ref{fact} gives \(l<p^*\). Moreover,
\eqref{eq:init-loop-endpoint-est} and~\eqref{eq:init-loop-apg-est} give
\[
l=f(\bm x_f)+\Delta,\qquad
\Delta>0,\qquad
\bm x_f,\bm x_g^{\rm ref}\in{\rm dom}(g_2),\qquad
g(\bm x_g^{\rm ref})-g^*\le\Delta\le\delta_{\rm init}\le1.
\]
By Line~\ref{alg:lb-compute-const}, the returned \(M_\Delta\) is the quantity
specified in~\eqref{eq:Mz-def}. Hence the returned tuple
satisfies Condition~\ref{cond:reference}.
The following theorem states the resulting \(\LBInit{}\) guarantee, including
the bounds in~\eqref{eq:MDelta-init-bound} and the operation count.

\begin{theorem}[\(\LBInit{}\) guarantee]
\label{thm:lb-total}

Suppose that Assumption~\ref{ass:basic} holds. Let \(D_g\) be
fixed as in~\eqref{eq:Dg-Omegag}, \(B_f \ge1\) be a valid
Lipschitz constant for \(f\) on \({\rm dom}(g_2)\), and
\(Q_{\rm init}\) be as in~\eqref{eq:m0-def}.
Then Algorithm~\ref{alg:lb-initialization} terminates and returns a tuple
satisfying Condition~\ref{cond:reference}, and the estimates
in~\eqref{eq:MDelta-init-bound} hold. Moreover, the total number of
basic operations is bounded by
\[
T_{\rm pre}:=
\cO\!\left(
\left[\sqrt{C_D L_{f_1}\Omega_g}\,Q_{\rm init}
+\sqrt{C_D L_{g_1}}\sqrt{Q_{\rm init}}\right]
{ \left(1+\log Q_{\rm init}\right)^2}
\right).
\]

\end{theorem}

\begin{proof}
It remains to prove
\eqref{eq:MDelta-init-bound} and the operation bound. 
Let \(j_0\) be the first index such that \(\Delta_{j_0}\le\theta\). Since the
stopping test is satisfied by iteration \(j_0\), the stopping index satisfies
\(K_0\le j_0\). Because \(\Delta_k=\delta_{\rm init}/2^k\) is decreasing, this
implies \(\Delta=\Delta_{K_0}\ge\Delta_{j_0}\). If \(j_0=0\), then
\(\Delta_{j_0}=\delta_{\rm init}\). If \(j_0\ge1\), the minimality of \(j_0\)
gives \(\Delta_{j_0-1}>\theta\), and hence
\(\Delta_{j_0}=\Delta_{j_0-1}/2>\theta/2\). Therefore
\[
\Delta^{-1}
\le
\max\left\{
\delta_{\rm init}^{-1},
\frac2\theta
\right\}
=
Q_{\rm init},
\]
where the equality follows from $\theta:=\min\{(p^*-f^*)/4,m_0/4\}$ and
the definition of \(Q_{\rm init}\) in~\eqref{eq:m0-def}. Thus \(\Delta\ge Q_{\rm init}^{-1}\). Furthermore, by
\eqref{eq:Dg-Omegag}, \(g(\bm x_f)-g^*+1\le\Omega_g\), while
\(g(\bm x_g^{\rm ref})\ge g^*\). Since \(\Omega_g\ge1\) and \(\Delta\le1\), equation \eqref{eq:MDelta-init-bound} follows because
\[
M_\Delta
=
\max\left\{
\frac{g(\bm x_f)-g(\bm x_g^{\rm ref})+1}{\Delta},1
\right\}
\le
\frac{\Omega_g}{\Delta}
\le
\Omega_g Q_{\rm init}.
\]

Finally, we count basic operations. Let \(T_k\) be the operation count of
iteration \(k\). The two APG calls in Lines~\ref{alg:lb-apgf}
and~\ref{alg:lb-apgg}, by Lemma~\ref{lem:conv-fistalinesearch}, require
\(\cO(D_g\sqrt{L_{f_1}/\Delta_k})\) and
\(\cO(D_g\sqrt{L_{g_1}/\Delta_k})\) operations, respectively. For the dual
call, the definition of \(M^{(k)}\) and~\eqref{eq:Dg-Omegag} give
\(M^{(k)}\le\Omega_g/\Delta_k\). We apply
Proposition~\ref{prop:conv-dual-algo} with
\[
\mu_{\phi_1}=\eta^{(k)}=\Delta_k/C_D,\qquad
M_z^\phi=M^{(k)},\qquad
\delta_{\rm D}^{(k)}=\Delta_k/(4C_DM^{(k)})
\]
and with \(L_{\phi_1}=L_{g_1}+\eta^{(k)}\), \(L_{\psi_1}=L_{f_1}\), and
\(\eta^{(k)}\le1\). The convention \(L_{g_1}\ge1\) absorbs the perturbation
term in \(L_{\phi_1}\). Under \(M^{(k)}\le\Omega_g/\Delta_k\), the
logarithmic factors in this proposition are
\({ \cO(1+\log(1/\Delta_k))}\), and the resulting dual-call cost is bounded by
the right-hand side in the estimate below. The two APG costs stated above are
bounded by the corresponding leading terms because \(C_D\ge D_g^2\),
\(\Omega_g\ge1\), and \(\Delta_k\le1\). Therefore
\[
T_k=
\cO\!\left(
\left[
\frac{\sqrt{C_D L_{f_1}\Omega_g}}{\Delta_k}
+
\frac{\sqrt{C_D L_{g_1}}}{\sqrt{\Delta_k}}
\right]
{ \left(1+\log\frac1{\Delta_k}\right)^2}
\right).
\]
Since \(\Delta_k=\delta_{\rm init}/2^k\ge\Delta\) for \(k\le K_0\) and
\(\Delta^{-1}\le Q_{\rm init}\),
\[
\sum_{k=0}^{K_0}\Delta_k^{-1}\le2\Delta^{-1}\le2Q_{\rm init},\qquad
\sum_{k=0}^{K_0}\Delta_k^{-1/2}
\le
\frac{\Delta^{-1/2}}{1-2^{-1/2}}
\le
\frac{\sqrt{Q_{\rm init}}}{1-2^{-1/2}}.
\]
Summing the bound on \(T_k\) over \(0\le k\le K_0\) gives the claimed
bound for \(T_{\rm pre}\).
\end{proof}

\section{Numerical Experiments}
\label{sec:experiment}

In this section, we apply \BiVFA{} and \SNtVFA{} to two SBO problems and
compare them with PB-APG \citep{chen2024penalty}, FBi-PG \citep{merchav2026dynamic}, IRE-APG \citep{shtern2026convergence}, and BiFPG \citep{bot2026accelerating}. For both experiments, we
set \(\epsilon=10^{-10}\) and \(\delta_{\rm init}=1\). The two proposed methods
use the same initialization routine and the prescribed adaptive inner tolerances. Their common
preprocessing is performed once, and its cost is included in the reported
runtime of each proposed method.

The experimental settings of the compared methods are described as follows.

\noindent\textup{(i)} PB-APG \citep{chen2024penalty} uses the penalty parameter \(10^4\).

\noindent\textup{(ii)} FBi-PG \citep{merchav2026dynamic} uses
\(\gamma=1.3\), with the regularization sequence
\(\alpha_k=(k+2)^{-1.3}\) for \(k\geq 0\).

\noindent\textup{(iii)} IRE-APG \citep{shtern2026convergence} uses
\(\beta=1\), corresponding to the regularization sequence
\(\sigma_k=k^{-1}\) for \(k\geq 1\).

\noindent\textup{(iv)} BiFPG \citep{bot2026accelerating} uses the authors' public
implementation\footnote{\url{https://github.com/echnen/fast-bilevel-methods}} with \(\alpha=4\), \(\sigma_e=\sigma_t=1\), \(c=100\),
\(\delta=1.9\), and the fixed stepsize
\(0.95/(L_{g_1}+\epsilon_0L_{f_1})\), where
\(\epsilon_0=100/2^{1.9}\).

All simulations are implemented using
MATLAB R2023a on a PC running Windows 11 with an AMD Ryzen 7 7840H CPU at
3.80 GHz and 16 GB RAM.

To illustrate the performance of the methods, we record the lower-level gap
and the absolute upper-level gap against runtime. The reported metrics are
\begin{equation}
\label{eq:experiment-metrics}
g(\bm x)-g^*,
\qquad
\lvert f(\bm x)-p^*\rvert.
\end{equation}
Here, for benchmarking purposes, we use the Greedy FISTA of \cite{liang2022improving} and the MATLAB function \texttt{fmincon} to
compute \(g^*\) and \(p^*\), respectively, to high accuracy.

\subsection{Linear Regression Problem with a Ball Constraint (LRPBC)}
\label{exper:lrpbc-l1}

In the first experiment, following \cite{merchav2023convex}, we consider a bilevel linear regression problem. The
lower-level problem fits the training data over an \(\ell_1\)-norm ball, while
the upper-level objective measures the validation loss. The problem takes the
following form:
\begin{equation}
\label{exper:lrpbc-l1-form}
\begin{array}{lcl}
&\min\limits_{\bm x\in\R^n}
& \frac12\|\bm A_{\rm val}\bm x-\bm b_{\rm val}\|^2\\
&{\rm s.t.}
& \bm x\in\argmin\limits_{\bm z\in C}
\left\{\frac12\|\bm A_{\rm tr}\bm z-\bm b_{\rm tr}\|^2\right\},
\end{array}
\end{equation}
where \(C=\{\bm x\in\R^n:\|\bm x\|_1\le1\}\), and
\((\bm A_{\rm tr},\bm b_{\rm tr})\) and
\((\bm A_{\rm val},\bm b_{\rm val})\) denote the training and validation data,
respectively. Related bilevel regression models can be found in
\citep{beck2014first,sabach2017first,jiang2023conditional,wang2024near}.

We conduct this experiment using the YearPredictionMSD data set\footnote{\url{https://archive.ics.uci.edu/dataset/203/yearpredictionmsd}}.
We randomly sample 1000 observations, and before dividing the data,
we add 50 dependent columns, each formed as a random linear combination of the original features. The data are then split into a training set containing 60\% of the samples and a validation set containing the remaining 40\%, and Gaussian noise with mean zero and standard deviation \(0.2\) is added to the validation data. The dependent columns make the training matrix rank deficient, which typically leads to multiple lower-level solutions.

\begin{figure}[htp]
\centering
\begin{subfigure}[b]{0.45\textwidth}
\centering
\includegraphics[width=\textwidth]{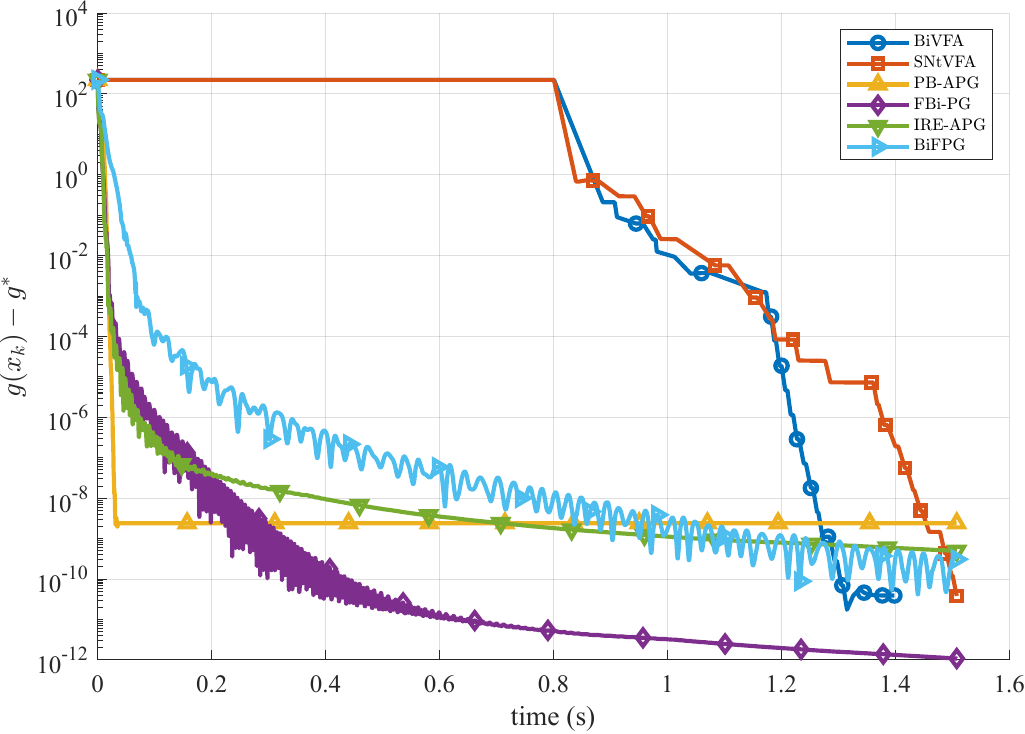}
\caption{Lower-level gap against runtime.}
\end{subfigure}
\hfill
\begin{subfigure}[b]{0.45\textwidth}
\centering
\includegraphics[width=\textwidth]{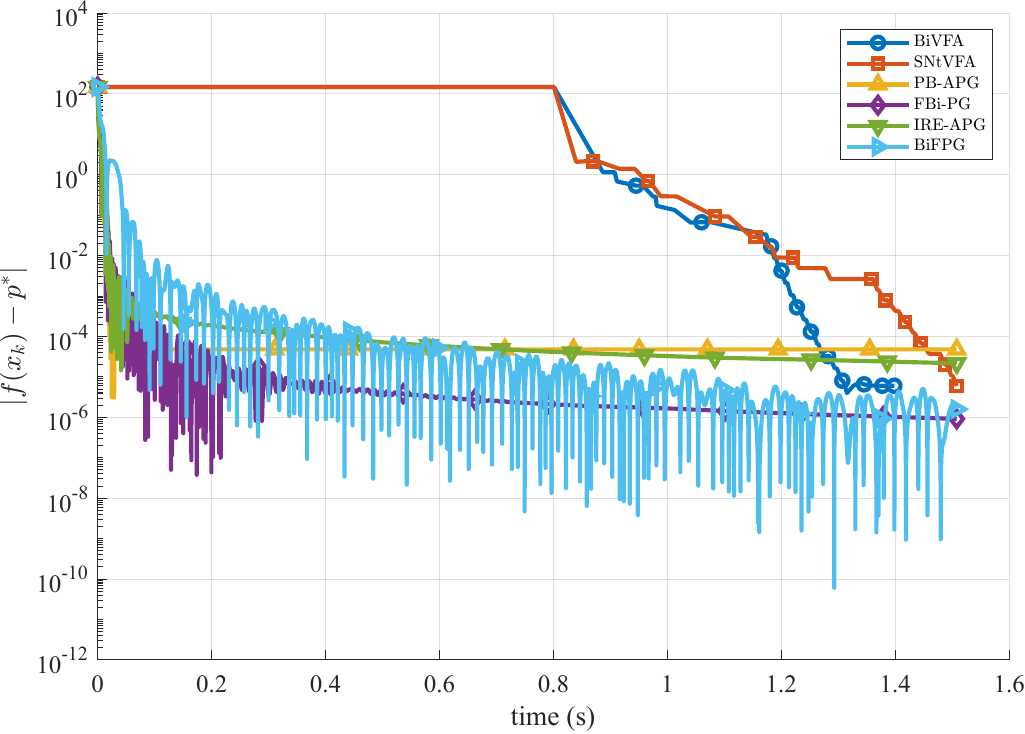}
\caption{Absolute upper-level gap against runtime.}
\end{subfigure}
\caption{Lower-level and absolute upper-level gaps over runtime for the LRPBC.}
\label{fig:lrpbc-l1-gaps}
\end{figure}

Figure~\ref{fig:lrpbc-l1-gaps} shows that some single-loop methods decrease
more rapidly in the low-accuracy regime. In the high-accuracy regime,
\BiVFA{} and \SNtVFA{} continue to reduce the reported residuals. Their
behavior in this regime is comparable with FBi-PG and remains competitive with
PB-APG, IRE-APG, and BiFPG under the reported parameter settings.
The initial flat segment corresponds to \(\LBInit{}\), which constructs the
reference certificate used by both methods. Accordingly, part of their reported
runtime reflects the cost of constructing this certificate and obtaining the
checkable oracle bounds that support the subsequent scalar decisions.
\FloatBarrier

\subsection{Integral Equation Problem (IEP)}
\label{exper:integral}

In the second experiment, following
\cite{sabach2017first,doron2023methodology}, we investigate an ill-conditioned
inverse problem arising from the discretization of a Fredholm integral
equation of the first kind~\citep{phillips1962technique}. The upper-level
objective selects a regularized solution from the lower-level least-squares
solution set. The problem takes the following form:
\begin{equation}
\label{exper:integral-form}
\begin{array}{lcl}
&\min\limits_{\bm x\in\R^n} & \bm x^{\T}\bm Q\bm x\\
&{\rm s.t.} & \bm x\in\argmin\limits_{\bm z\in C}
\left\{\frac12\|\bm A\bm z-\bm b\|^2\right\},
\end{array}
\end{equation}
where \(C=\R^n_+\) and \(\bm Q=\bm L^{\T}\bm L+\bm I\), with
\(\bm L=\texttt{get\_l}(150,1)\). The matrix \(\bm A_{\rm base}\), the
noiseless vector \(\bm b_{\rm true}\), and the true solution are generated by
\texttt{phillips(100)} from the Regularization Tools package\footnote{\url{http://www2.imm.dtu.dk/~pcha/Regutools/}}.
Since \(C=\R^n_+\) is unbounded, this instance falls outside the
compact-domain setting of Assumption~\ref{ass:basic}. We nevertheless include
it to examine the empirical behavior of the proposed methods beyond that
setting.
We set \(\bm b=\bm b_{\rm true}+0.2\bm w\), where \(\bm w\) is sampled from
the standard normal distribution. To make the lower-level problem rank
deficient, we form \(\bm A=\bm A_{\rm base}[\bm I,\bm D]\), where \(\bm D\)
deterministically selects 50 evenly spaced Phillips columns. This construction
makes the lower-level least-squares objective non-strictly convex and typically
leads to multiple lower-level solutions.

\begin{figure}[htp]
\centering
\begin{subfigure}[b]{0.45\textwidth}
\centering
\includegraphics[width=\textwidth]{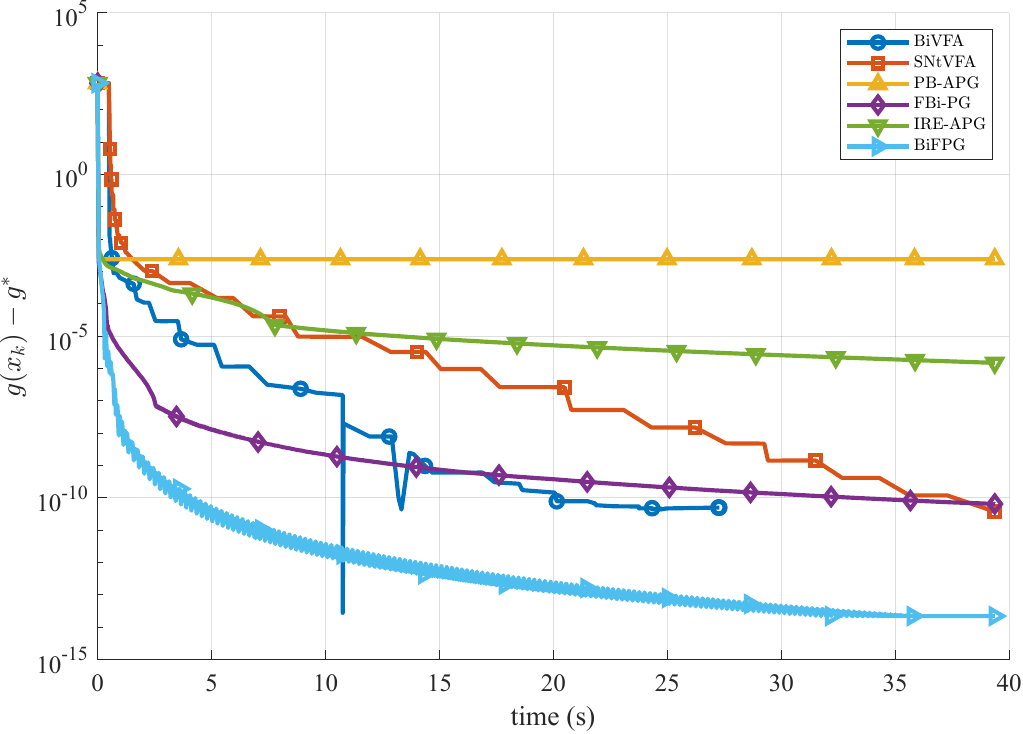}
\caption{Lower-level gap against runtime.}
\end{subfigure}
\hfill
\begin{subfigure}[b]{0.45\textwidth}
\centering
\includegraphics[width=\textwidth]{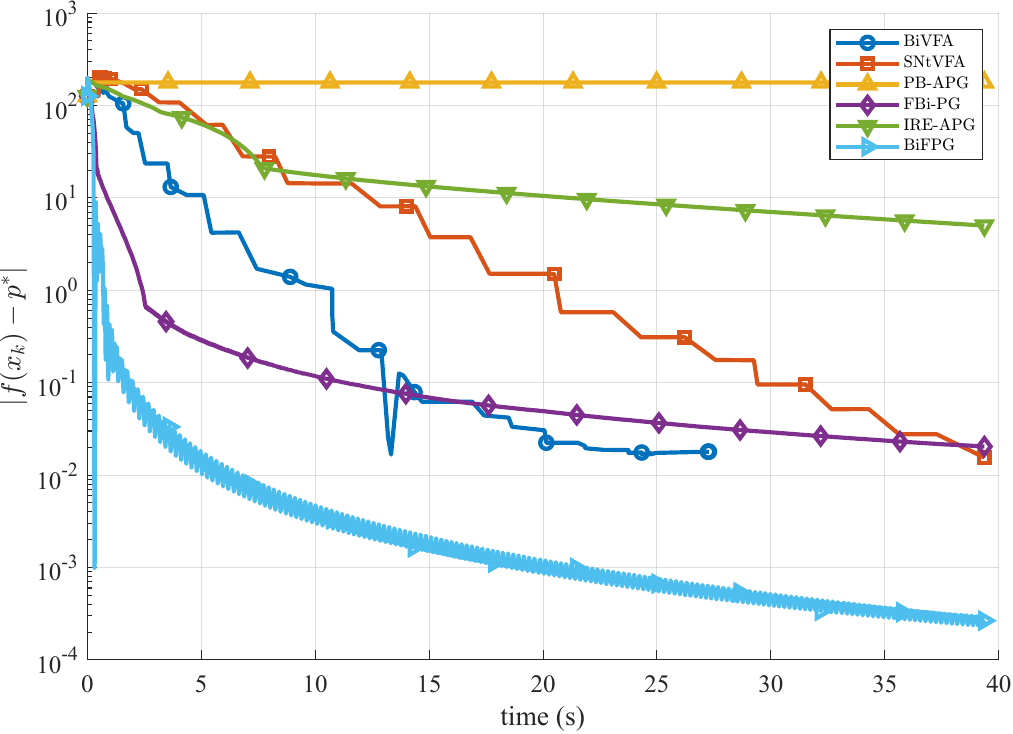}
\caption{Absolute upper-level gap against runtime.}
\end{subfigure}
\caption{Lower-level and absolute upper-level gaps over runtime for the IEP.}
\label{fig:iep-gaps}
\end{figure}

Figure~\ref{fig:iep-gaps} shows that several single-loop methods reduce the
reported gaps more rapidly at early runtimes. The reported runtimes of our
methods include the one-time \LBInit{} cost. The associated initial flat segment
is less pronounced than in Figure~\ref{fig:lrpbc-l1-gaps} because the IEP runs
extend over a longer time horizon. \BiVFA{} and \SNtVFA{} eventually attain
lower-level gaps of order \(10^{-10}\), whereas their upper-level gaps remain of
order \(10^{-2}\) on this instance. Under the reported parameter settings, their
final lower- and upper-level gaps are comparable in order of magnitude to those
of FBi-PG and
smaller than those of PB-APG and IRE-APG. BiFPG attains the smallest final gaps
on this instance.

\section{Conclusion}
\label{sec:conclusion}

This paper developed a value-function root-finding approach with checkable oracle bounds
for composite convex simple bilevel optimization. By constructing a
value-function oracle, the proposed framework replaces level-set proximal access
with feasibility and value-error bounds together with multiplier information under a
proximal-composite operation model. The resulting bisection and safeguarded
Newton-type algorithms compute an \((\epsilon,\epsilon)\)-optimal solution.
After a one-time reference-tuple initialization, both
scalar-search phases have \(\tO(\epsilon^{-1/2})\) accuracy dependence. The
numerical results support the predicted scalar-search behavior on
representative composite convex instances.

\begin{appendices}

\section{\texorpdfstring{Generic Dual Oracle}{Generic Dual Oracle}}
\label{app:dual-oracle-details}

This appendix proves the generic dual-oracle guarantee stated in
Proposition~\ref{prop:conv-dual-algo}. Throughout,
Assumption~\ref{ass:dual-subproblem} is in force for
Problem~\eqref{p:func-cons}, and we use the notation introduced in
Subsection~\ref{sec:generic-subroutines}.

\subsection{\texorpdfstring{Basic Estimates for the Dual Function}{Basic Estimates for the Dual Function}}

The compact domain gives the following modulus.
\begin{lemma}
\label{lem:lip-psi}
There exists a constant \(B_\psi{\ge1}\) such that \(\psi\) is
\(B_\psi\)-Lipschitz continuous on \({\rm dom}(\phi_2)\).
\end{lemma}
\begin{proof}
Compactness of \({\rm dom}(\phi_2)\) and continuity of \(\nabla\psi_1\) give
\(B_{\psi_1}:=\max_{\bm x\in{\rm dom}(\phi_2)}\|\nabla\psi_1(\bm x)\|<\infty\).
The set \({\rm dom}(\phi_2)\) is convex, so the mean-value inequality gives
\({ \lvert\psi_1(\bm x)-\psi_1(\bm y)\rvert}\le B_{\psi_1}\|\bm x-\bm y\|\) for all
\(\bm x,\bm y\in{\rm dom}(\phi_2)\). Combining this with the assumed
\(l_{\psi_2}\)-Lipschitz continuity of \(\psi_2\) on \({\rm dom}(\phi_2)\)
shows that {\(B_\psi:=\max\{1,B_{\psi_1}+l_{\psi_2}\}\)} is valid.
\end{proof}
Below, \(B_\psi\) denotes any known Lipschitz constant satisfying
\(B_\psi\ge1\).
The next lemma adapts the monotonicity and Lipschitz estimates of
\citep[Lemma~3.2]{xu2022first} to the present composite Lagrangian
\(\cL_\phi\) and its parametric minimizer \(\bm x_\phi(z)\).
\begin{lemma}
\label{lem:lips-xzphi}
For any \(z_1,z_2 \in [0, 2M_z^\phi)\), we have
\begin{equation}
\label{eq:mono-phi}
(z_1-z_2)\bigl(\psi(\bm x_\phi(z_1))-\psi(\bm x_\phi(z_2))\bigr)
\le -\mu_{\phi_1}\|\bm x_\phi(z_1)-\bm x_\phi(z_2)\|^2.
\end{equation}
\begin{equation}
\label{eq:lip-xzphi}
\|\bm x_\phi(z_1)-\bm x_\phi(z_2)\| \le \frac{B_\psi}{\mu_{\phi_1}}\lvert z_1-z_2\rvert.
\end{equation}
\end{lemma}
\begin{proof}
For \(i=1,2\), write \(\bm x_i=\bm x_\phi(z_i)\). The \(\mu_{\phi_1}\)-strong
convexity of \(\cL_\phi(\cdot,z_i)\), defined in~\eqref{p:dzphi}, gives
\[
\phi(\bm x_1)+z_1\psi(\bm x_1)
\le
\phi(\bm x_2)+z_1\psi(\bm x_2)-\frac{\mu_{\phi_1}}2\|\bm x_1-\bm x_2\|^2,
\]
\[
\phi(\bm x_2)+z_2\psi(\bm x_2)
\le
\phi(\bm x_1)+z_2\psi(\bm x_1)-\frac{\mu_{\phi_1}}2\|\bm x_1-\bm x_2\|^2.
\]
Adding these inequalities gives~\eqref{eq:mono-phi}. By
Lemma~\ref{lem:lip-psi},
\[
\mu_{\phi_1}\|\bm x_1-\bm x_2\|^2
\le
{ \lvert z_1-z_2\rvert\,\lvert\psi(\bm x_1)-\psi(\bm x_2)\rvert}
\le
B_\psi { \lvert z_1-z_2\rvert}\,\|\bm x_1-\bm x_2\|.
\]
If \(\bm x_1=\bm x_2\), then~\eqref{eq:lip-xzphi} is immediate; otherwise,
dividing by \(\|\bm x_1-\bm x_2\|\) gives~\eqref{eq:lip-xzphi}.
\end{proof}

The next lemma records the Danskin derivative formula for the dual function,
including the right derivative at the endpoint \(z=0\).
\begin{lemma}
\label{lem:derivative-dz}
The dual function \(d_\phi(z)\), defined in~\eqref{p:dzphi}, is differentiable
on \((0,2M_z^\phi)\) and satisfies
\begin{equation}\label{eq:main-derivative}
d_\phi'(z)=\psi(\bm x_\phi(z)),
\qquad \forall z\in(0,2M_z^\phi).
\end{equation}
Moreover, the right derivative at \(z=0\) exists and satisfies
\(d_{\phi,+}'(0)=\psi(\bm x_\phi(0))\).
\end{lemma}
\begin{proof}
For \(z\in(0,2M_z^\phi)\), the Danskin formula established in the proof of
\cite[Proposition~A.22]{bertsekas1971control}, applied to the minimization
envelope \(d_\phi\), gives \eqref{eq:main-derivative}.

For the endpoint \(z=0\), let \(\bm x_t:=\bm x_\phi(t)\) for
\(t>0\) and \(\bm x_0:=\bm x_\phi(0)\). Since \(\bm x_t\) minimizes
\(\cL_\phi(\cdot,t)\), we have
\(
d_\phi(t)
=
\phi(\bm x_t)+t\psi(\bm x_t)\). 
Since \(\bm x_0\) minimizes \(\cL_\phi(\cdot,0)=\phi(\cdot)\), we have
\(d_\phi(0)=\phi(\bm x_0)\) and
\(\phi(\bm x_t)\ge \phi(\bm x_0)=d_\phi(0)\). Hence
\[
d_\phi(t)-d_\phi(0)
=
\phi(\bm x_t)-d_\phi(0)+t\psi(\bm x_t)
\ge t\psi(\bm x_t).
\]
On the other hand, using the comparison point \(\bm x_0\) in the
minimization defining \(d_\phi(t)\) gives
\[
d_\phi(t)
\le
\phi(\bm x_0)+t\psi(\bm x_0)
=
d_\phi(0)+t\psi(\bm x_0).
\]
Combining the last two displays and dividing by \(t>0\), we obtain
\[
\psi(\bm x_t)
\le
\frac{d_\phi(t)-d_\phi(0)}{t}
\le
\psi(\bm x_0).
\]
The Lipschitz estimate~\eqref{eq:lip-xzphi}, applied with \(z_1=t\) and
\(z_2=0\), gives
\(\|\bm x_t-\bm x_0\|\le (B_\psi/\mu_{\phi_1})t\). Hence
\(\bm x_t\to \bm x_0\) as \(t\downarrow0\).
Since \(\psi\) is continuous on \({\rm dom}(\phi_2)\),
letting \(t\downarrow0\) gives
\(
d'_{\phi,+}(0)=\psi(\bm x_\phi(0)).
\)
\end{proof}

The stationarity guarantee determines the sign of the exact dual derivative.
\begin{lemma}
\label{lem:sign-dzphi}
Let \(\tilde z\in[0,2M_z^\phi)\) and \(0<\delta\le 1/B_\psi\). Set
\(\varepsilon_1=\mu_{\phi_1}\delta^2\) and
\(\varepsilon_2=\delta\). If \(\tilde{\bm x}\in{\rm dom}(\phi_2)\) satisfies
\({\rm dist}(\bm 0,\partial_{\bm x}\cL_{\phi}(\tilde{\bm x},\tilde z))\le\varepsilon_1\),
then

\begin{equation}
\label{eq:muphi}
\mu_{\phi_1}\|\tilde{\bm x} -\bm x_\phi(\tilde z)\| \le \varepsilon_1.
\end{equation}
Moreover, if \([\psi(\tilde{\bm x})]_+\le\varepsilon_2\), then
\([\psi(\bm x_\phi(\tilde z))]_+\le2\varepsilon_2\). Otherwise,
\(\psi(\bm x_\phi(\tilde z))>0\).
\end{lemma}
\begin{proof}

Choose \(\xi\in\partial_{\bm x}\cL_{\phi}(\tilde{\bm x},\tilde z)\) with
\(\|\xi\|\le\varepsilon_1\). With
\(\bm 0 \in \partial_{\bm x} \cL_{\phi}(\bm x_\phi(\tilde z),\tilde z)\),
strong monotonicity of the subdifferential of the
\(\mu_{\phi_1}\)-strongly convex function
\(\cL_\phi(\cdot,\tilde z)\)
\citep[Theorem~5.24(iii)]{beck2017first} gives
\[
\mu_{\phi_1}\|\tilde{\bm x} -\bm x_\phi(\tilde z)\|
\le \|\xi\| \le \varepsilon_1,
\]
which proves~\eqref{eq:muphi}.
The Lipschitz continuity of \(\psi\) on \({\rm dom}(\phi_2)\) therefore gives
\begin{equation}
\label{eq:err-psi-hatx}
\lvert\psi(\tilde{\bm x})-\psi(\bm x_\phi(\tilde z))\rvert
\le B_\psi \|\tilde{\bm x}-\bm x_\phi(\tilde z)\|
\overset{\eqref{eq:muphi}}{\le}
\frac{B_\psi}{\mu_{\phi_1}}\varepsilon_1
=B_\psi\delta^2
\le\delta
=\varepsilon_2,
\end{equation}
where the last inequality follows from \(\delta\le1/B_\psi\). Since
\([\cdot]_+\) is \(1\)-Lipschitz, \eqref{eq:err-psi-hatx} gives
\([\psi(\bm x_\phi(\tilde z))]_+\le2\varepsilon_2\) whenever
\([\psi(\tilde{\bm x})]_+\le\varepsilon_2\). If
\([\psi(\tilde{\bm x})]_+>\varepsilon_2\), the same estimate gives
\(\psi(\bm x_\phi(\tilde z))>0\).
\end{proof}

\subsection{\texorpdfstring{Interval Search}{Interval Search}}

For fixed \(\delta\), we say that \(Z=[a,b]\subset[0,2M_z^\phi)\) satisfies
the endpoint conditions if there exist
\(\tilde{\bm x}_a,\tilde{\bm x}_b\in{\rm dom}(\phi_2)\) such that
\begin{equation}
\label{eq:endpoint-conditions}
\begin{gathered}
{\rm dist}(\bm 0,\partial_{\bm x}\cL_\phi(\tilde{\bm x}_a,a))\le\varepsilon_1,
~
[\psi(\tilde{\bm x}_a)]_+>\varepsilon_2,\quad\text{ and }\quad
{\rm dist}(\bm 0,\partial_{\bm x}\cL_\phi(\tilde{\bm x}_b,b))\le\varepsilon_1,
~
[\psi(\tilde{\bm x}_b)]_+\le\varepsilon_2.
\end{gathered}
\end{equation}

With this terminology, Lemmas~\ref{lem:sign-dzphi}
and~\ref{lem:derivative-dz} lead to the interval-search routine
Algorithm~\ref{alg:inter-searchphi} (\IntV{}). The interval-search routine adapts the one-dimensional multiplier search of
\citep[Algorithm~3 and Lemma~3.8]{xu2022first} to the present composite
dual-oracle setting. 
The routine either returns a KKT
pair for~\eqref{p:func-cons} or a multiplier interval satisfying the endpoint
conditions, and Lemma~\ref{lem:output-interphi} below records its output
guarantee and operation bound.

\begin{algorithm}[ht!]
\caption[Interval search (\IntV{})]{Interval search:
\(({\rm flag},\tilde{\bm x},\tilde z,Z,\tilde{\bm x}_a,\tilde{\bm x}_b)\)
\(=\IntV(\phi,\psi,\sigma,M_z^\phi,\bm x_{\rm init};\delta)\)}
\label{alg:inter-searchphi}
\DontPrintSemicolon
\KwIn{\((\phi,\psi)\) satisfies Assumption~\ref{ass:dual-subproblem},
\(0<\sigma\le M_z^\phi\), \(0<\delta\le1/B_\psi\), and
\(\bm x_{\rm init}\in{\rm dom}(\phi_2)\).}
Set \(\varepsilon_1,\varepsilon_2,\varepsilon_3\) as in~\eqref{eq:dual-tolerances}.\;
Invoke \(\tilde{\bm x}{\gets}
\APGm(\phi_1,\phi_2,\varepsilon_1,\bm x_{\rm init})\).\;
\If{\([\psi(\tilde{\bm x})]_{+}\le \varepsilon_2\)}{
\Return \({\rm flag}=1\), \((\tilde{\bm x},0)\).
\label{alg:inter-searchphi-returnz-zero}
\Comment*[r]{zero-multiplier KKT pair}
}
{ Set} \(a{\gets}0\), \(b{\gets}\sigma\), and \(\tilde{\bm x}_a{\gets}\tilde{\bm x}\).\;
Invoke \(\tilde{\bm x}_b{\gets}
\APGm
(\phi_1 + b\psi_1,\phi_2 + b\psi_2,\varepsilon_1,\bm x_{\rm init})\).\;
\While{\([\psi(\tilde{\bm x}_b)]_{+}>\varepsilon_2\) and \(b < M_z^\phi\)}{\label{alg:inter-searchphi-step7s}
{ Set} \(a{\gets} b\), \(b{\gets}2b\), and \(\tilde{\bm x}_a{\gets}\tilde{\bm x}_b\).\;
Invoke \(\tilde{\bm x}_b{\gets}
\APGm
(\phi_1 + b\psi_1,\phi_2 + b\psi_2,\varepsilon_1,\bm x_{\rm init})\).\;
}
\If{\(\lvert b \psi(\tilde{\bm x}_b)\rvert\le \varepsilon_3\)}{
\Return \({\rm flag}=1\), \((\tilde{\bm x}_b,b)\).
\label{alg:inter-searchphi-returnz-b}
\Comment*[r]{positive-multiplier KKT pair}
}
\Else{
\Return \({\rm flag}=0\), \(Z=[a,b]\), \((\tilde{\bm x}_a,\tilde{\bm x}_b)\).
\label{alg:inter-searchphi-returnz-ab}
\Comment*[r]{endpoint conditions for bisection}
}
\end{algorithm}

\begin{lemma}
\label{lem:output-interphi}
Suppose that Assumption~\ref{ass:dual-subproblem} holds for
Problem~\eqref{p:func-cons}.
After at most \(T_{\rm int}\) basic operations,
Algorithm~\ref{alg:inter-searchphi} produces either an
\((\varepsilon_1,\varepsilon_2,\varepsilon_3)\)-KKT pair
\((\tilde{\bm x},\tilde z)\) of Problem~\eqref{p:func-cons} with
\(\tilde z\in [0,2M_z^\phi)\), or an interval
\(Z=[a,b]\subset[0,2M_z^\phi)\) that contains an optimal multiplier of
Problem~\eqref{p:func-cons}, where
\(\varepsilon_1,\varepsilon_2,\varepsilon_3\) are defined
by~\eqref{eq:dual-tolerances}, and \(T_{\rm int}\) satisfies
\begin{equation}
\label{eq:barT}
T_{\rm int}
=
\cO\left(
\left(
1+
\left\lceil\log_2 \frac{M_z^\phi}{\sigma}\right\rceil_+
\right)
\sqrt{\frac{L_{\phi_1}+M_z^\phi L_{\psi_1}}{\mu_{\phi_1}}}
\left(
1+
\left\lceil
\log_2\frac{(L_{\phi_1}+M_z^\phi L_{\psi_1})D_{\phi}}{\mu_{\phi_1}\delta^2}
\right\rceil_+
\right)
\right).
\end{equation}
In the second case, the returned interval and endpoint points
\(\tilde{\bm x}_a,\tilde{\bm x}_b
\in{\rm dom}(\phi_2)\) satisfy the endpoint
conditions in~\eqref{eq:endpoint-conditions}.
\end{lemma}
\begin{proof}
\emph{Complexity.} Each \(\APGm\) call uses tolerance \(\varepsilon_1\) and a multiplier
\(z\in[0,2M_z^\phi)\). Hence the smooth part in that call has Lipschitz
constant at most \(L_{\phi_1}+2M_z^\phi L_{\psi_1}\). The initial call at
\(z=0\), the call at \(z=\sigma\), and the doubling loop make at most
\(2+\lceil\log_2(M_z^\phi/\sigma)\rceil_+\) calls in total. The APG guarantee in
Lemma~\ref{lem:apgmu-complexity} then gives~\eqref{eq:barT}.

\emph{Correctness.}
Let \(\bar z_\phi\le M_z^\phi\) be a dual solution from
Assumption~\ref{ass:dual-subproblem}, and set
\(\chi(z):=\psi(\bm x_\phi(z))\). Since \(d_\phi\) is concave and
\(M_z^\phi\) lies to the right of a maximizer,
Lemma~\ref{lem:derivative-dz} gives \(\chi(M_z^\phi)\le0\). If
Line~\ref{alg:inter-searchphi-returnz-zero} is reached, the APG stopping rule
at \(z=0\), the branch condition
\([\psi(\tilde{\bm x})]_+\le\varepsilon_2\), and
\({ \lvert0\cdot\psi(\tilde{\bm x})\rvert}=0\) give an
\((\varepsilon_1,\varepsilon_2,\varepsilon_3)\)-KKT pair. 

If the positive
multiplier return in Line~\ref{alg:inter-searchphi-returnz-b} is reached, APG
gives stationarity at \((\tilde{\bm x}_b,b)\). At termination of the doubling
loop in Line~\ref{alg:inter-searchphi-step7s}, we have either
\([\psi(\tilde{\bm x}_b)]_+\le\varepsilon_2\), or \(b\ge M_z^\phi\). In the
second case, \eqref{eq:mono-phi} shows that \(\chi\) is nonincreasing, and
therefore \(\chi(b)\le\chi(M_z^\phi)\le0\). Hence
the APG stopping rule at \(b\), together with~\eqref{eq:muphi}, gives
\(\|\tilde{\bm x}_b-\bm x_\phi(b)\|\le\varepsilon_1/\mu_{\phi_1}\). Since
\(\psi\) is \(B_\psi\)-Lipschitz on \({\rm dom}(\phi_2)\),
\[
\psi(\tilde{\bm x}_b)
\le
\psi(\bm x_\phi(b))
+B_\psi\|\tilde{\bm x}_b-\bm x_\phi(b)\|
=
\chi(b)+B_\psi\|\tilde{\bm x}_b-\bm x_\phi(b)\|
\le
\frac{B_\psi\varepsilon_1}{\mu_{\phi_1}}.
\]
Since \(\varepsilon_1=\mu_{\phi_1}\delta^2\), \(\varepsilon_2=\delta\), and
\(\delta\le1/B_\psi\), the last display gives
\([\psi(\tilde{\bm x}_b)]_+\le\varepsilon_2\).
The complementarity condition
\({ \lvert b\psi(\tilde{\bm x}_b)\rvert}\le\varepsilon_3\) in
Line~\ref{alg:inter-searchphi-returnz-b} gives the required KKT pair.

It remains to consider the interval return. The APG stopping rule and endpoint
updates give the endpoint conditions~\eqref{eq:endpoint-conditions}. At
\(b\), these conditions, \eqref{eq:muphi}, and the Lipschitz
continuity of \(\psi\) give
\begin{equation}
\label{eq:b-endpoint-bounds}
\begin{gathered}
{\rm dist}(\bm 0,\partial_{\bm x}\cL_\phi(\tilde{\bm x}_b,b))\le\varepsilon_1,
\qquad
[\psi(\tilde{\bm x}_b)]_+\le\varepsilon_2=\delta,\\
{ \lvert\psi(\tilde{\bm x}_b)-\chi(b)\rvert}
\le B_\psi\|\tilde{\bm x}_b-\bm x_\phi(b)\|
\overset{\eqref{eq:muphi}}{\le}
B_\psi\varepsilon_1/\mu_{\phi_1}
\le\delta.
\end{gathered}
\end{equation}
Here the last inequality uses
\(\varepsilon_1=\mu_{\phi_1}\delta^2\) and \(\delta\le1/B_\psi\). For the
lower endpoint, \eqref{eq:endpoint-conditions}, \eqref{eq:muphi}, and
Lipschitz continuity yield
\[
\chi(a)=\psi(\bm x_\phi(a))\ge
\psi(\tilde{\bm x}_a)-B_\psi\|\tilde{\bm x}_a-\bm x_\phi(a)\|
> \delta-B_\psi\delta^2\ge0,
\]
so \(\chi(a)>0\).
Assume, for contradiction, that \([a,b]\) contains no optimal multiplier. By
\eqref{eq:mono-phi}, \(\chi\) is nonincreasing, and
Lemma~\ref{lem:lip-psi} and \eqref{eq:lip-xzphi} imply that \(\chi\) is
continuous. Moreover, by concavity of \(d_\phi\) and
Lemma~\ref{lem:derivative-dz}, any zero of \(\chi\) in \((a,b]\) is a dual
maximizer. Since \(\chi(a)>0\), the contradiction hypothesis gives
\(\chi(b)>0\). As
\(\chi(M_z^\phi)\le0\) and \(\chi\) is nonincreasing, we also have
\(b<M_z^\phi\). The bound in~\eqref{eq:b-endpoint-bounds} then gives
\(-\delta<\psi(\tilde{\bm x}_b)\le\delta\). 
Since \(b<M_z^\phi\),
\[
{ \lvert b\psi(\tilde{\bm x}_b)\rvert}\le M_z^\phi\delta
\le 4M_z^\phi\delta=\varepsilon_3.
\]
Together with the first two bounds in
\eqref{eq:b-endpoint-bounds}, Line~\ref{alg:inter-searchphi-returnz-b}
would return the positive-multiplier KKT pair, not an interval, a contradiction. Therefore the returned
interval contains an optimal multiplier.
\end{proof}

\subsection{\texorpdfstring{Multiplier Bisection}{Multiplier Bisection}}

When Algorithm~\ref{alg:inter-searchphi} (\IntV{}) returns an interval
\(Z=[a_{\rm init},b_{\rm init}]\), together with endpoint points
\(\tilde{\bm x}_a,\tilde{\bm x}_b
\in{\rm dom}(\phi_2)\), Algorithm~\ref{alg:bisec-multiphi}
(\Bisec{}) adapts the single-constraint multiplier-bisection idea of
\citep[Algorithm~4]{xu2022first} to refine this interval in the present
composite dual-oracle setting while preserving the endpoint
conditions~\eqref{eq:endpoint-conditions} from
Lemma~\ref{lem:output-interphi}. The next lemma records the resulting KKT
residual guarantee and operation bound.

\begin{algorithm}[ht]
\caption[Multiplier bisection (\Bisec{})]{Multiplier bisection:
\((\tilde{\bm x},\tilde z)=\Bisec
(\phi,\psi,Z,M_z^\phi,\tilde{\bm x}_a^0,\tilde{\bm x}_b^0,\bm x_{\rm init};\delta)\)}
\label{alg:bisec-multiphi}
\DontPrintSemicolon
\KwIn{\(Z=[a_{\rm init},b_{\rm init}]\) satisfies the endpoint conditions
in~\eqref{eq:endpoint-conditions}, with endpoint points
\(\tilde{\bm x}_a^0,\tilde{\bm x}_b^0\) returned by \IntV{},
\(0<\delta\le1/B_\psi\), and
\(\bm x_{\rm init}\in{\rm dom}(\phi_2)\).}
Set \(\varepsilon_1,\varepsilon_2,\varepsilon_3\) as in~\eqref{eq:dual-tolerances}, and initialize
\(a{\gets} a_{{\rm init}}\), \(b{\gets} b_{{\rm init}}\),
\(\tilde{\bm x}_a{\gets}\tilde{\bm x}_a^0\), and
\(\tilde{\bm x}_b{\gets}\tilde{\bm x}_b^0\).\;
\While{\(b-a>\mu_{\phi_1}\delta^3\)}{\label{alg:bisec-multiphi-while}
{ Set} \(z_{\rm mid}{\gets}(a+b)/2\).\;
Invoke \(\tilde{\bm x}_{\rm mid}{\gets}
\APGm
(\phi_1 + z_{\rm mid}\psi_1,\phi_2 + z_{\rm mid}\psi_2,
\varepsilon_1,\bm x_{\rm init})\).\;
\If{\([\psi(\tilde{\bm x}_{\rm mid})]_+>\varepsilon_2\)}{
{ Set} \(a{\gets} z_{\rm mid}\) and \(\tilde{\bm x}_a{\gets}\tilde{\bm x}_{\rm mid}\).\;
}
\ElseIf{\([\psi(\tilde{\bm x}_{\rm mid})]_+ \le \varepsilon_2\) and
\(\lvert z_{\rm mid} \psi(\tilde{\bm x}_{\rm mid})\rvert\le \varepsilon_3\)}{
\Return \((\tilde{\bm x},\tilde z)=(\tilde{\bm x}_{\rm mid}, z_{\rm mid})\).
\label{alg:bisec-multiphi-step6}
\Comment*[r]{KKT pair}
}
\Else{
{ Set} \(b{\gets} z_{\rm mid}\) and \(\tilde{\bm x}_b{\gets}\tilde{\bm x}_{\rm mid}\).\label{alg:bisec-multiphi-step11}
}
}
\Return \((\tilde{\bm x},\tilde z)=(\tilde{\bm x}_b,b)\).
\label{alg:bisec-multiphi-step13}
\Comment*[r]{terminal endpoint return}
\end{algorithm}

\begin{lemma}
\label{lem:epsiKKT-output0phi}
Suppose that Assumption~\ref{ass:dual-subproblem} holds for
Problem~\eqref{p:func-cons}.
Given an interval
\(Z=[a_{\rm init},b_{\rm init}]\subset[0,2M_z^\phi)\)
satisfying the endpoint conditions in~\eqref{eq:endpoint-conditions}, with the
corresponding endpoint points
\(\tilde{\bm x}_a^0,\tilde{\bm x}_b^0\) supplied to the algorithm,
Algorithm~\ref{alg:bisec-multiphi} produces an
\((\varepsilon_1,\varepsilon_2,\varepsilon_3)\)-KKT pair
\((\tilde{\bm x},\tilde z)\) of
Problem~\eqref{p:func-cons} with
\(\tilde z\in[0,2M_z^\phi)\) after at most \(T_{\rm bis}\) basic
operations, where
\begin{equation}\label{eq:barTbisec}
T_{\rm bis} \le
\cO\left(\sqrt{\frac{L_{\phi_1}+M_z^\phi L_{\psi_1}}{\mu_{\phi_1}}}
\left(
1+
\left\lceil
\log_2\frac{(L_{\phi_1}+M_z^\phi L_{\psi_1})D_{\phi}}{\mu_{\phi_1}\delta^2}
\right\rceil_+
\right)
\left(
1+
\left\lceil
\log_2\frac{M_z^\phi}{\mu_{\phi_1}\delta^3}
\right\rceil_+
\right)\right).
\end{equation}
\end{lemma}
\begin{proof}
The two endpoint-update branches preserve
\eqref{eq:endpoint-conditions}: an infeasible midpoint becomes the lower
endpoint, while a feasible midpoint that fails the complementarity test becomes
the upper endpoint. If Line~\ref{alg:bisec-multiphi-step6} returns a midpoint,
the APG stationarity guarantee and the branch tests already give
\eqref{eq:kkt-tol}. The returned multiplier also lies in
\([0,2M_z^\phi)\), since the midpoint lies in the current interval.

The stopping rule and \(b_{\rm init}-a_{\rm init}<2M_z^\phi\) give
\(\cO(1+\lceil\log_2(M_z^\phi/(\mu_{\phi_1}\delta^3))\rceil_+)\) bisection
iterations. For each midpoint,
\(L_{\phi_1}+z_{\rm mid}L_{\psi_1}\le L_{\phi_1}+2M_z^\phi L_{\psi_1}\), so
Lemma~\ref{lem:apgmu-complexity} gives~\eqref{eq:barTbisec}, after absorbing
numerical constants in \(\cO(\cdot)\).

It remains to verify the terminal endpoint return. If no midpoint is returned,
Algorithm~\ref{alg:bisec-multiphi} returns \((\tilde{\bm x}_b,b)\) with
\(b\in[0,2M_z^\phi)\) and
\(b-a\le\mu_{\phi_1}\delta^3\). By the maintained endpoint
conditions~\eqref{eq:endpoint-conditions} and \(\varepsilon_2=\delta\)
from~\eqref{eq:dual-tolerances}, both endpoints are
\(\varepsilon_1\)-stationary, with \(\psi(\tilde{\bm x}_a)>\delta\) and
\(\psi(\tilde{\bm x}_b)\le\delta\). For \(z=a,b\),
Lemma~\ref{lem:sign-dzphi} at \(z\), Lipschitz continuity of \(\psi\),
\(\varepsilon_1=\mu_{\phi_1}\delta^2\), and \(\delta\le1/B_\psi\) give
\({ \lvert\psi(\tilde{\bm x}_z)-\psi(\bm x_\phi(z))\rvert}\le\delta\). Combining this
estimate with Lemma~\ref{lem:lips-xzphi} yields
\[
{ \lvert\psi(\tilde{\bm x}_b)-\psi(\tilde{\bm x}_a)\rvert}
\le\delta+B_\psi\|\bm x_\phi(b)-\bm x_\phi(a)\|+\delta
\le2\delta+\frac{B_\psi^2}{\mu_{\phi_1}}(b-a)
\le2\delta+B_\psi^2\delta^3
\le3\delta.
\]
Together with \(\psi(\tilde{\bm x}_a)>\delta\) and
\(\psi(\tilde{\bm x}_b)\le\delta\), this gives
\(
-2\delta\le\psi(\tilde{\bm x}_b)\le\delta.
\)
Hence
\([\psi(\tilde{\bm x}_b)]_+\le\varepsilon_2\) and
\({ \lvert\psi(\tilde{\bm x}_b)\rvert}\le2\delta\). Since \(b<2M_z^\phi\),
\(
{ \lvert b\psi(\tilde{\bm x}_b)\rvert}\le4M_z^\phi\delta=\varepsilon_3.
\)
Together with the \(\varepsilon_1\)-stationarity bound
\({\rm dist}(\bm 0,\partial_{\bm x}\cL_\phi(\tilde{\bm x}_b,b))
\le\varepsilon_1\), these bounds give~\eqref{eq:kkt-tol} for that pair.
\end{proof}

We now combine the interval-search and multiplier-bisection guarantees to prove
Proposition~\ref{prop:conv-dual-algo}.

\begin{proof}[Proof of Proposition~\ref{prop:conv-dual-algo}]
If \(\IntV\) returns \({\rm flag}=1\),
Lemma~\ref{lem:output-interphi} gives an
\((\varepsilon_1,\varepsilon_2,\varepsilon_3)\)-KKT pair with
\(\tilde z\in[0,2M_z^\phi)\). If \({\rm flag}=0\),
Lemma~\ref{lem:output-interphi} returns an interval satisfying the endpoint
conditions required by Lemma~\ref{lem:epsiKKT-output0phi}; hence \(\Bisec\) returns the
same KKT residual bounds and multiplier range. Adding the operation bounds in
Lemmas~\ref{lem:output-interphi} and~\ref{lem:epsiKKT-output0phi} gives the
displayed operation bound. 
\end{proof}

\end{appendices}


\begingroup
\bibliographystyle{abbrv}
\bibliography{ref_sbo.bib}
\endgroup

\end{document}